\documentclass{amsart}
\usepackage{amssymb,amsmath}
\usepackage[american,french]{babel}
\usepackage{amsopn}
\usepackage{fancyhdr}
\title{Ergodicity for  the stochastic Complex Ginzburg--Landau equations}
\author{Cyril ODASSO \\ \\  Ecole Normale Sup\'erieure de Cachan, antenne de Bretagne,\\ Avenue Robert Schuman,
 Campus de Ker Lann, 35170 Bruz (FRANCE). \\ and
\\ IRMAR,  UMR 6625 du CNRS, Campus de Beaulieu,  35042 Rennes cedex (FRANCE)}

\newtheorem{Theorem}{Theorem}[section]

\newtheorem{Proposition}[Theorem]{Proposition}
\newtheorem{Lemma}[Theorem]{Lemma}
\newtheorem{Corollary}[Theorem]{Corollary}
\newtheorem{Remark}[Theorem]{Remark}

\makeatletter
\newif\ifmsbmloaded@

\@addtoreset{equation}{section}

\makeatother

\def\R{\mathbb R}
\def\N{\mathbb N}

\def\C{\mathbb C}
\def\E{\mathbb E}
\def\P{\mathbb P}
\def\Pcal{\mathcal{P}}
\def\Hcal{\mathcal{H}}

\def\Dr{\mathcal D}
\def\Br{\mathcal B}

\def\F{\mathcal F}

\newcommand{\BLANC}[1]{   }

\newcommand{\abs}[1]{\left\vert#1\right\vert}
\newcommand{\norm}[1]{\left\Vert#1\right\Vert}
\newcommand{\eps}{\varepsilon}
\newcommand{\sig}{\sigma}
\renewcommand{\i}{\textrm{i}}
\def \Espace{\renewcommand{\arraystretch}{1.7} }
\newcommand{\carre}{    \begin{flushright}
                $\Box$
            \end{flushright}}

\newcommand{\abssig}[1]{\abs{#1}_{2\sig +2}^{2\sig+2}}
\newcommand{\absquatre}[1]{\abs{#1}_{4\sig +2}^{4\sig+2}}
\newcommand{\B}[1]{ \ensuremath{\abs{#1}^{2\sig}  #1}}

\newcommand{\chal}[1]{ \ensuremath{ \frac{d#1}{dt} + ( \eps + \i )A#1}}

\newcommand{\chalcarreH}[2]{ \ensuremath{ \frac{d\norm{#1}^2}{dt} %
 + #2 \eps \norm{#1}_2^2}}
\newcommand{\glc}[1]{ \ensuremath{ \chal{#1}+( \eta +\lambda \i ) \B{#1} }}

\newcommand{\glcDelta}[1]{\ensuremath{\frac{d#1}{dt}-(\eps+\i)\Delta#1 %
+( \eta + \lambda\i ) \B{#1} }}

\begin{document}

\setcounter{page}{84}

\selectlanguage{american}
\maketitle
\pagestyle{fancy}





\noindent\textbf{Abstract}:
 We study a stochastic complex Ginzburg--Landau (CGL) equation driven by a smooth noise in space
 and we establish exponential convergence of the Markov transition semi-group toward a unique invariant
 probability measure. Since Doob Theorem does not seem not to be useful in our situation, a coupling method is used.
  In order to make
this method easier to understand, we first focus on two simple examples which contain most of the arguments  and the
 essential difficulties.

\

\noindent\textbf{R\'esum\'e}:
 Nous consid\'erons l'\'equation de Ginzburg--Landau Complexe bruit\'ee par un bruit blanc en temps et r\'egulier par
rapport aux variables spatiales et nous \'etablissons le caract\`ere exponentiellement m\'elangeant
 du semi-groupe de Markov vers une unique mesure de probabilit\'e invariante.
Comme le Th\'eor\`eme de Doob semble ne pas pouvoir \^etre appliquer, nous utilisons une m\'ethode dite de couplage.
 Pour une meilleur compr\'ehension, nous focaliserons d'abord notre attention sur deux exemples qui bien que tr\`es simples
contiennent l'essentiel des difficult\'es.

\


\


\section*{Introduction}

Originally introduced to describe a phase transition in superconductivity \cite{GL}, the Complex Ginzburg--Landau (CGL)
 equation also models the propagation of dispersive non-linear waves in various areas of physics such as hydrodynamics
\cite{Newel1}, \cite{Newel2}, optics, plasma physics, chemical reaction \cite{Huber}...

 When working in non-homogenous
or random media, a noise is often introduced and the stochastic CGL equation may be more representative than the
deterministic one.

The CGL equation arises in the same areas of physics as the non-linear Schr\"odinger (NLS) equation. In fact, the CGL
equation  is obtained by adding two viscous terms to the NLS equation. The inviscid limits of the deterministic and
stochastic CGL equation to the NLS equation are established in \cite{Bebouche} and \cite{KSlimit}, respectively. The
 stochastic NLS equation is studied in \cite{dBD1} and \cite{dBD2}.

Ergodicity of the stochastic CGL equation is established in \cite{Marc} when the noise is invertible and in
\cite{H} for the one-dimensionnal cubic case when the noise is  diagonal, does not
 depend on the solution and is smooth in space.

Our aim in this article is to study ergodicity for stochastic CGL equation under very general assumptions.

Let us recall that the stochastic CGL equation has the form
$$
\left \{
\begin{array}{rcll}
\glcDelta{u} & = & b(u)\frac{dW}{dt} ,& \\
u(t,x) & = & 0,           &\mbox{   for   } x \in \delta D, \; t>0 \\
u(0,x) & = & u_0(x), &\mbox{   for   } x \in D.
\end{array}
\right .
$$
The unknown $u$ is a complex valued process depending on $x \in D$, $D \subset \R^d$ a bounded domain, and $t\geq 0$.

We want to consider noises which may be degenerate and our work is in the spirit of  \cite{BKL}, \cite{EMS}, \cite{H}, \cite{KS00}, \cite{KS},
 \cite{KS2}, \cite{KS3},  \cite{Matt} and \cite{S}. Many ideas of this article are taken from
 these works. However, we develop several generalisations.

The main idea is to compensate the degeneracy of the noise on some subspaces by dissipativity arguments, the so-called
Foias-Prodi estimates. A coupling method is developped in a sufficiently general framework to be applied and prove
exponential convergence to equilibrium.

To describe the ideas, it is convenient to introduce $(e_k)_{k\in \N^*}$ the eigenbasis of the operator $-\Delta$ with
Dirichlet boundary conditions (if periodic boundary conditions were considered, it would be the Fourier basis) and $P_N$
 the eigenprojector on the first $N$ modes.

The main assumption of the papers cited above as well as in this work is that the noise is non-degenerate on the space
 spanned by $(e_k)_{1\leq k\leq N}$ for $N$ sufficiently large. In \cite{H}, \cite{KS3} and \cite{S},
 the noise is also additive, i.e. $b(u)$ does not depend on $u$. The method developped in \cite{Matt}
allows to treat more general noises and, in \cite{Matt}, $b$ is allowed to depend on $P_N u$. However, in this latter work,
 the author restricts his attention to the case when the high modes are not perturbed by noise. It is claimed that
 the method can be generalized to treat a noise which hits all components.
Such a generalisation is contained in \cite{Matt2} in the purely additive case.

Here we develop also such a generalization and treat a noise which may hit all modes but depends only on $P_Nu$.
We have chosen to use ideas both from \cite{Matt} and from \cite{KS3}, \cite{S}. We hope that this makes our proof easier
to understand. Moreover, we get rid of the assumption that $b$ is diagonal in the basis $(e_k)_{k\in \N^*}$.

Also, if we work in the space $L^2(D)$, it is not difficult to get a Lyapunov structure and Foias-Prodi estimates. Thus,
 with an additive noise or with a noise as in \cite{Matt}, our results would be a rather easy applications of these
methods.

However, this works only for small values of $\sig$, namely $\sig < \frac{2}{d}$. It is well known that the CGL equations
 are also well-posed for $\sig \in \left[ \frac{2}{d},\frac{2}{d-2} \right)$
 ( $\sig \in \left[\frac{2}{d},\infty \right)$  for $ d \in \{1,2\}$) provided we work with
$H^1(D)$--valued solutions and the nonlinearity is defocusing $(\lambda=1)$. We also develop the coupling method in
 that context and show that it is possible to find a convenient Lyapunov structure and derive Foias-Prodi estimates.
Thus we prove exponential convergence to equilibrium for the noises described above in all the cases when it is known that
there exists a unique global solution and an invariant measure.

Moreover, using the smoothing effect of CGL and an interpolation argument, we are able to prove exponential convergence in
the Wasserstein norm in $H^s(D)$ for any $s<2$. This give convergence to equilibrium for less regular funtionnal.

In order to make the understanding of the method easier, we start with two simple examples which motivate and introduce
all arguments in a simpler context. The first example is particulary simple. It introduces the idea of coupling and
the use of Girsanov transform to construct a coupling. The second example is similar to the one considered in  \cite{Matt}.
  However, it contains further difficulties and more details are given. We have tried to isolate every key argument.
This is also the opportunity to state a very general result giving conditions implying exponential mixing
(Theorem \ref{Th_Theo}). It is a strong generalization of Theorem 3.1 of \cite{KS3}.

Then, in section 2 we deal with CGL equations. We state and prove the general ergodicity result described above


\section{Preliminary results}

The proof of our result is obtained by the combination of two main ideas: the coupling and the Foias-Prodi estimate. The first subsection is a simple example devoted to understand the use of the notion of coupling.
The second subsection is a two dimensional example devoted to understand how we use the two main ideas. The third
subsection is the statement of an abstract result which is both fundamental and technical. The other subsections are devoted to the proof of
this abstract result. The understanding of the proof of the abstract result is not necessary to the understanding of the
rest of the article. On the contrary the three first subsections contain the main ideas of this article.

\subsection{A simple example}
\

In this subsection we introduce the notion of coupling and we motivate it on a simple example.

Let $\Pi$ the one-dimensionnal torus. We consider the following example. We
denote by $X(., x_0)$ the unique solution in $\Pi$ of
\begin{equation}\label{Eq_coup_ex}
\frac{dX}{dt}+f(X)=\frac{dW}{dt}, \quad X(0, x_0)=x_0,
\end{equation}
where $f:\Pi \to \R $ is a Lipschitz function and $W$ is a one-dimensionnal brownian motion. It is easy to prove that
 $X$ is a Markovian process. We denote by $(\Pcal_t)_t$ its Markovian transition semigroup.

We recall the definition of $\norm{\mu}_{var}$, the total variation of a
finite real measure $\mu$:
 $$
\norm{\mu}_{var}= \sup  \left\{\abs{\mu(\Gamma)}\;|\; \Gamma \in \Br(\Pi) \right\},
$$
where we denote by $\Br(\Pi)$ the set of the Borelian subsets of $\Pi$. It is well known that $\norm{.}_{var}$ is the dual
 norm of $\abs{.}_\infty$.
We prove that there exists a unique invariant measure $\nu$ and that for any probability measure $\mu$
$$
\norm{\Pcal_t^* \mu -\nu}_{var}\leq c e^{-\beta t}.
$$
Using a completeness argument and the markovian property of $X$, we obtain that it is sufficient to prove
that for any $\psi:\Pi\to \R$ borelian bounded and for any $(t,x_1, x_2) \in \R^+\times\Pi^2$, we have
$$
\abs{\E\psi(X(t,x_1))-\E\psi(X(t,x_2))}\leq c\abs{\psi}_\infty e^{-\beta t}.
$$
Clearly it is sufficient to find $(X_1(t),X_2(t))$  such that
for any $(i,t)\in \{1,2\}\times\R^+$, we have $\Dr(X_i(t))=\Dr(X(t,x_i))$, where $\Dr$ means
 distribution,  and
\begin{equation}\label{Eq_coup_ex_c}
\abs{\E\psi(X_1(t))-\E\psi(X_2(t))}\leq c\abs{\psi}_\infty e^{-\beta t}.
\end{equation}

Now we introduce the notion of coupling.
Let $(\mu_1,\mu_2)$ be two distributions on a same space $(E,\mathcal{E})$. Let $(\Omega,\mathcal{F},\P)$ be a probability
 space and let $(Z_1,Z_2)$ be two random variables $(\Omega,\mathcal{F}) \to (E,\mathcal{E})$. We say that $(Z_1,Z_2)$
is a coupling of  $(\mu_1,\mu_2)$ if $\mu_i=\Dr(Z_i)$ for $i=1,2$.

\begin{Remark}
Although the marginal laws of $(Z_1,Z_2)$ are imposed, we have a lot of freedom when choosing the law of the couple $(Z_1,Z_2)$.
 For instance, let us consider  $(W_1,W_2)$ a two-dimensional brownian motion. Let $\mu$ be the Wiener measure on $\R$, which
 means that $\mu=\Dr(W_1)=\Dr(W_2)$. Then $(W_1,W_2)$, $(W'_1,W'_2)=(W_1,W_1)$ and $(W''_1,W''_2)=(W_1,-W_1)$ are
 three couplings of $(\mu,\mu)$. These three couplings have very different laws. In the one hand, $W_1$ and $W_2$ are
 independent and $W_1 \not = \pm W_2 \textrm{ a.s.} $ and in the other hand $W'_1=W'_2$ and $W''_1=-W''_2$.
\end{Remark}

In order to establish \eqref{Eq_coup_ex_c}, we remark that it is sufficient to build
$\left(X_1,X_2\right)$ a coupling of $\left(\Dr\left(X(\cdot,x_1)\right),\Dr\left(X(\cdot,x_2)\right)\right)$ on $\R^+$
 such that for any $t\geq 0$
\begin{equation}\label{Eq_coup_ex_e}
\P\left(X_1(t)\not = X_2(t)\right)\leq c e^{-\beta t}.
\end{equation}
By induction, it suffices to construct a coupling on a fixed interval $[0,T]$. Indeed, we first set
$$
X_i(0)=x_i \; , \; i=1,2.
$$
Then we build a probability space $\left( \Omega', \F', \P' \right)$ and a measurable function
 $(\omega',t,x_1,x_2)\to Z_i(t,x_1,x_2)$ such that for any $(x_1,x_2)$, $(Z_i(\cdot,x_1,x_2))_{i=1,2}$ is a coupling of
 $(X(\cdot,x_i))_{i=1,2}$ on $[0,T]$.

The induction argument is then as follows. Assuming that we have built $(X_1,X_2)$ on $[0,nT]$, we take $(Z_1,Z_2)$ as
 above independant of $(X_1,X_2)$ on $[0,nT]$ and set
$$
X_i(nT+t)=Z_i(t,X_1(nT),X_2(nT)), \quad \textrm{ for } t \in (0,T].
$$
The Markov property of $X$ implies that $(X_1,X_2)$ is a coupling of $(\Dr(X(\cdot,x_1)),\Dr(X(\cdot,x_2)))$ on $[0,(n+1)T]$.

The coupling $(Z_1,Z_2)$ on $[0,T]$ constructed below satisfies the following properties
\begin{eqnarray}
\P\left(Z_1(T,x_1,x_2)= Z_2(T,x_1,x_2)\right)&\geq p_0 >0, & \textrm{ if } x_1 \not = x_2, \label{Eq_coup_ex_f}\\
\P\left(Z_1(.,x_1,x_2)= Z_2(.,x_1,x_2)\right)&=1, & \textrm{ if } x_1  = x_2.\label{Eq_coup_ex_g}
\end{eqnarray}
Invoking \eqref{Eq_coup_ex_g}, we obtain that
$$
\P\left(X_1(nT)\not=X_2(nT)| X_1((n-1)T)=X_2((n-1)T) \right)=0.
$$
Thus it follows
$$
\Espace
\begin{array}{lcl}
\P\left( X_1(nT)\not=X_2(nT) \right) & \leq & \P\left(  X_1((n-1)T)\not=X_2((n-1)T)\right)\times \\
&&\P\left(X_1(nT)\not=X_2(nT)| X_1((n-1)T)\not=X_2((n-1)T) \right) .
\end{array}
$$
We easily get from \eqref{Eq_coup_ex_f} and  \eqref{Eq_coup_ex_g}
$$
\P\left(  X_1(t)\not = X_2(t),\; \textrm{ for some } t\geq nT \right)\leq (1-p_0)^n,
$$
which implies \eqref{Eq_coup_ex_e} and allows us to conclude.

Before building $(Z_1,Z_2)$ such that \eqref{Eq_coup_ex_f} and \eqref{Eq_coup_ex_g} hold, we need to
define some notions.
 Let $\mu$, $\mu_1$ and $\mu_2$ be three probability measures on a space $(E,\mathcal{E})$ such that $\mu_1$ and $\mu_2$ are absolutely continuous
 with respect to $\mu$. We set
$$
\Espace
\begin{array}{lcl}
d\abs{\mu_1 - \mu_2}&=&\abs{\frac{d\mu_1}{d\mu}-\frac{d\mu_2}{d\mu}} d\mu, \\
d(\mu_1 \wedge \mu_2)&=&(\frac{d\mu_1}{d\mu}\wedge\frac{d\mu_2}{d\mu}) d\mu,\\
d(\mu_1 - \mu_2)^+&=&(\frac{d\mu_1}{d\mu}-\frac{d\mu_2}{d\mu})^+ d\mu.
\end{array}
$$
These definitions do not depend on the choice of $\mu$.
Moreover we have
\begin{equation}\label{Eq_coup_ex_b}
\norm{\mu_1-\mu_2}_{var}= \frac{1}{2} \abs{\mu_1 - \mu_2}(E)= (\mu_1 - \mu_2)^+(E)=\frac{1}{2} \int_E  \abs{ \frac{d\mu_1}{d\mu}-\frac{d\mu_2}{d\mu}}d\mu.
\end{equation}
The following Lemma is the key of our proof.
\begin{Lemma}\label{lem_norm_var_sens}
Let $(\mu_1,\mu_2)$ be two probability measures on $(E,\mathcal{E})$. Then
$$
\norm{\mu_1-\mu_2}_{var}= \min \P(Z_1\not = Z_2).
$$
The minimum is taken over the coupling $(Z_1,Z_2)$ of $(\mu_1,\mu_2)$. Such a coupling exists and is called a maximal coupling
and has the following property:
$$
\P(Z_1=Z_2, Z_1 \in \Gamma)=(\mu_1 \wedge \mu_2)(\Gamma)\; \textrm{ for any } \Gamma \in \mathcal{E}.
$$
\end{Lemma}
The proof of Lemma \ref{lem_norm_var_sens} is given in the Appendix.
We consider $W'$ a Wiener process. If $x_1=x_2=x$, we choose the trivial coupling $(Z_i(.,x,x))_{i=1,2}$ on $[0,T]$.
In other words, we set $Z_1(.,x,x)=Z_2(.,x,x)=X'(.,x)$ on $[0,T]$ where $X'(.,x)$ is the solution of \eqref{Eq_coup_ex}
associated with $W'$. Thus \eqref{Eq_coup_ex_g} is clear.

For $x_1\not =x_2$,
the idea is borrowed from \cite{KS3}. We consider
($\widetilde{Z}_1(.,x_1,x_2)$,\linebreak \noindent $Z_2(.,x_1,x_2)$) the maximal coupling of
 $\left(\Dr(X(\cdot,x_1) +\frac{T-\cdot}{T}(x_2-x_1)),\Dr(X(.,x_2))\right)$ on $[0,T]$ and we set
$Z_1(t,x_1,x_2)=\widetilde{Z}_1(t,x_1,x_2) -\frac{T-t}{T}(x_2-x_1)$. Then it is easy to see that
 $(Z_i(.,x_1,x_2))_{i=1,2}$ is a coupling of $(\Dr(X(.,x_i))_{i=1,2}$ on $[0,T]$ and we have
\begin{equation}\label{Eq_coup_ex_i}
\P\left(Z_1(T,x_1,x_2) = Z_2(T,x_1,x_2) \right)\geq \P\left(\widetilde{Z}_1(.,x_1,x_2)= Z_2(.,x_1,x_2) \right).
\end{equation}
We need the following result which is lemma D.1 of \cite{Matt}
\begin{Lemma}\label{lem_tech_coup_inf}
Let $\mu_1$ and $\mu_2$ be two probability measures on a space $(E,\mathcal{E})$. Let $A$ be an event of $E$. Assume that
 $\mu_1^A=\mu_1(A\cap.)$ is equivalent to $\mu_2^A=\mu_2(A\cap.)$.  Then for any $p>1$ and $C>1$
 $$
\int_A \left ( \frac{d \mu_1^A}{d \mu_2^A} \right )^{p+1} d\mu_2 \leq C <\infty
\quad \textrm{implies} \quad
\left(\mu_1\wedge\mu_2\right)(A) \geq  \left( 1-\frac{1}{p} \right)\left( \frac{\mu_1(A)^p}{pC} \right)^{\frac{1}{p-1}}.
$$
\end{Lemma}
Using \eqref{Eq_coup_ex_i} and Lemma \ref{lem_norm_var_sens} and \ref{lem_tech_coup_inf} with $E=C([0,T];\Pi)$, we obtain that
\begin{equation}\label{Eq_coup_ex_o}
\P\left(Z_1(T,x_1,x_2) = Z_2(T,x_1,x_2) \right)\geq \left( 1-\frac{1}{p} \right)
\left( p\int_E \left ( \frac{d \widetilde{\mu}_1}{d \mu_2} \right )^{p+1} d\mu_2  \right)^{-\frac{1}{p-1}},
\end{equation}
where $(\widetilde{\mu}_1,\mu_2)=(\Dr(X(\cdot,x_1) +\frac{T-\cdot}{T}(x_2-x_1)),\Dr(X(.,x_2)))$ on $[0,T]$.

We use a Girsanov formula to estimate $\int_E \left ( \frac{d \widetilde{\mu}_1}{d \mu_2} \right )^{p+1} d\mu_2 $.
Setting $\widetilde{X}(t)=X(t,x_1)+\frac{T-t}{T}(x_2-x_1)$, we obtain that $\widetilde{\mu}_1$ is the distribution of
$\widetilde{X}$ under the probability
$\P$ and that $\widetilde{X}$ is the unique solution of
$$
\frac{d\widetilde{X}}{dt} -\frac{1}{T}(x_2-x_1)+f(\widetilde{X}(t)+\frac{T-t}{T}(x_2-x_1))=\frac{dW}{dt}, \quad \widetilde{X}(0)=x_2.
$$
We set $W'(t)=W(t)+\int_0^t d(s) dt$, where
\begin{equation}\label{Eq_un_sept_bis}
d(t)=\frac{1}{T}(x_2-x_1)+f(\widetilde{X}(t))-f(\widetilde{X}(t)+\frac{T-t}{T}(x_2-x_1)).
\end{equation}
Then $\widetilde{X}$ is a solution of
\begin{equation}\label{Eq_coup_ex_l}
\frac{d\widetilde{X}}{dt}+f(\widetilde{X})=\frac{dW'}{dt}, \quad \widetilde{X}(0)=x_2,
\end{equation}
We are working on the torus and $f$ is continuous, therefore $d$ is uniformly bounded:
$$\abs{d(t)}\leq \frac{1}{T}+2\abs{f}_\infty.$$
Hence, the Novikov condition is satisfied and the Girsanov formula can be applied. Then we set
$$
d\P'= \exp\left(\int_0^t d(s) dW(s)-\frac{1}{2} \int_0^t \abs{d(s)}^2 dt  \right)d\P
$$
 We deduce from the Girsanov formula that $\P'$ is a probability measure under which $W'$ is a
 brownian motion and $\widetilde{X}$ is a solution of \eqref{Eq_coup_ex_l}, then the law of $\widetilde{X}$ under $\P'$
 is $\mu_2$. Moreover
\begin{equation}\label{Eq_coup_ex_lbis}
\int_E \left ( \frac{d \widetilde{\mu}_1}{d \mu_2} \right )^{p+1} d\mu_2 \leq
\exp\left( c_p \left( \frac{1}{T}+\abs{f}_\infty^2 T \right) \right),
\end{equation}
which allows us to conclude this example. Indeed, by applying  \eqref{Eq_coup_ex_i}, \eqref{Eq_coup_ex_o} and \eqref{Eq_coup_ex_lbis}
we get \eqref{Eq_coup_ex_f}.

\subsection{A representative two-dimensionnal example}
\

The example we consider now is a two dimensional system which mimics the decomposition of a stochastic partial differential
 equation according to low and high modes of the solution. This example allows the introduction of the main ideas in a
simplified context, the system has the form
\begin{equation}\label{Eq_un_neuf}
\left \{
\Espace \begin{array}{lcl}
d X+ 2 X dt +f(X,Y)dt & = & \sig_l(X)d\beta , \\
d Y+ 2 Y dt +g(X,Y)dt & = &    \sig_h(X)d\eta  , \\
\lefteqn{ X(0)=x_0 , \quad Y(0)=y_0.}
\end{array}
\right .
\end{equation}
We set $u=\left( X, Y \right)$ and $W= \left( \beta , \eta \right)$.
We use the following assumptions
\begin{equation}\label{Eq_un_neuf_bis}
 \left \{
\Espace \begin{array}{ll}
\textrm{i)} & f, \; g, \; \sig_l \textrm{ and } \sig_h \textrm{ are bounded and Lipschitz}, \\
\textrm{ii)}& \textrm{There exists } K_0>0 \textrm{ such that,}  \\
& f(x,y)x+g(x,y)y \geq -(\abs{x}^2+\abs{y}^2+K_0),\quad (x,y)\in \R^2 .
\end{array}
\right.
\end{equation}
Condition i) ensures existence and uniqueness of a solution to \eqref{Eq_un_neuf} once the initial data $u_0=(x_0,y_0)$ is given.
 It is also classical that weak existence and uniqueness holds. We denote by $X(\cdot,u_0)$, $Y(\cdot,u_0)$, $u(\cdot,u_0)$
the solution where $u_0=(x_0,y_0)$ and $u=(X,Y)$. Moreover, it is easy to see that, by ii), there exists an invariant measure $\nu$.

Contrary to section 1.1, we want to allow degenerate noises. More precisely, we want to treat the case when the noise on
 the second equation may vanish. This possible degeneracy is compensated by a dissipativity assumption. We use the
 following assumptions.
\begin{equation}\label{Eq_un_neuf_ter}
 \left \{
\Espace \begin{array}{ll}
\textrm{i)} & \textrm{There exists } \sig_0>0 \textrm{ such that,} \; \sig_l(x)\geq \sig_0, \; x\in \R.\\
\textrm{ii)}& \abs{g(x,y_1)-g(x,y_2)}\leq \abs{y_1-y_2},\quad (x,y_1,y_2)\in \R^2.
\end{array}
\right.
\end{equation}
By the dissipativity method (see \cite{DPZ1} section 11.5), ii) implies exponential convergence to equilibrium for the
second equation if $X$ is fixed. Whilst the coupling argument explained in section 1.1 can be used to treat the first
 equation when $Y$ is fixed.
Note however that we need a more sophisticated coupling here. Indeed, the simple coupling explained above seems to be usefull
 only for additive noise.

Here, we explain how these two arguments may be coupled to treat system \eqref{Eq_un_neuf}.
The essential tool which allows to treat system \eqref{Eq_un_neuf} is the so-called Foias-Prodi estimate which reflects the
dissipativity property of the second equation. It is a simple consequence of \eqref{Eq_un_neuf_ter}ii)
\begin{Proposition}\label{Prop_Foais_Ex}
Let $(u_i,W_i)_{i=1,2}$ be two weak solutions of \eqref{Eq_un_neuf} such that
$$
X_1(s)=X_2(s), \quad \eta_1(s)=\eta_2(s), \quad s\in [0,t],
$$
then
$$
\abs{u_1(t)-u_2(t)}\leq \abs{u_1(0)-u_2(0)}e^{-t}
$$
\end{Proposition}
Since the noise on the second equation might be degenerate, there is no hope to use Girsanov formula on the full system. We can
use it to modify the drift of the first equation only and it is not possible to derive a strong estimate as
\eqref{Eq_coup_ex_e}.

Recall that in section 1.1, we have built the coupling $(X_1,X_2)$ of $(\Dr(X(\cdot,x_0^1)),\Dr(X(\cdot,x_0^2)))$ by induction on
$[0,kT]$ by using a coupling $(Z_i(\cdot,x_0^1,x_0^2))_{i=1,2}$ of $(\Dr(X(\cdot,x_0^i)))_{i=1,2}$ on $[0,T]$ which
satisfies \eqref{Eq_coup_ex_g}. Then if $(X_1,X_2)$ were coupled at time $kT$, $(X_1,X_2)$ would be coupled
 on $[kT,\infty)$ with probability one. Thus to conclude, it was sufficient to establish  \eqref{Eq_coup_ex_f}.

In this section, since we couple $(X_1,X_2)$, but not $(Y_1,Y_2)$, then there is no hope that a couple $(X_1,X_2)$ coupled
 at time $kT$ remains coupled at time $(k+1)T$ with probability one.

However, coupling the $X$'s and using Foias-Prodi estimates, we obtain a coupling $(u_1,u_2)$ of
 $(\Dr(u(\cdot,u_0^1)), \Dr(u(\cdot,u_0^2)))$ on $\R^+$ such that
\begin{equation}\label{Eq_etoile}
\P\left(\abs{u_1(t)-u_2(t)}>c e^{-\beta t}  \right)\leq c e^{-\beta t}(1+\abs{u_0^1}^2+\abs{u_0^2}^2).
\end{equation}
This estimate does not imply the decay of the total variation of $\Pcal_t^*\delta_{u_0^1}-\Pcal_t^*\delta_{u_0^2}$, but
 the decay of this quantity in the Wasserstein distance $\abs{\cdot}_{Lip_b}^*$ which is the dual norm of the lipschitz
 and bounded functions. Indeed, for $\psi$ lipschitz and bounded, we clearly have
$$
\Espace \begin{array}{rcl}
\abs{\E \psi(u(t,u_0^1))-\E \psi(u(t,u_0^2))} & = & \abs{\E \psi(u_1(t))-\E \psi(u_2(t))} , \\
        & \leq & 2 \abs{\psi}_\infty \P\left(\abs{u_1(t)-u_2(t)}>c e^{-\beta t}  \right)+\abs{\psi}_{Lip} c e^{-\beta t},
\end{array}
$$
and then by \eqref{Eq_etoile}
\begin{equation}\label{Eq1_15bis}
\abs{\E \psi(u(t,u_0^1))-\E \psi(u(t,u_0^1))} \leq c \abs{\psi}_{Lip_b} e^{-\beta t}(1+\abs{u_0^1}^2+\abs{u_0^2}^2).
\end{equation}
The idea of the proof is the following. We couple $(\Dr(X(\cdot,u_0^i),\eta))_{i=1,2}$. Then using the Foias-Prodi estimate,
 we control $Y_1-Y_2$ which is equivalent to control $u_1-u_2$. By controlling $u_1-u_2$, we control the
probability to remain coupled.
\begin{Remark}
In the general case $f$, $g$  are not globally lipschitz and bounded and a cut-off has to be used.
 This further difficulty will be treated in the context of the CGL equation below.
\end{Remark}
It is convenient to introduce the following functions:
$$
l_0(k)=\min\left\{l \in \{0,...,k\} | P_{l,k}\right\},
$$
 where $\min \phi =\infty$ and
$$
(P_{l,k})
\left \{
\Espace
\begin{array}{l}
X_1(t)=X_2(t), \quad \eta_1(t)=\eta_2(t), \quad \forall\; t \in [lT,kT],\\
\abs{u_i(lT)}\leq d^* ,\quad i=1,2.
\end{array}
\right.
$$
The first requirement in $(P_{l,k})$ states that the two solutions of the first equation are coupled on $[lT,kT]$.
Notice that Proposition \ref{Prop_Foais_Ex} gives
\begin{equation}\label{Eq_abis}
l_0(k)=l \quad \textrm{implies} \quad   \abs{u_1(t)-u_2(t)}\leq 2d^* e^{-(t-lT)}, \; \textrm{ for any } t \in [lT,kT].
\end{equation}
 From now on we say that $(X_1,X_2)$ are coupled at $kT$ if $l_0(k)\leq k$, in other words
 if $l_0(k)\not = \infty$.

We set
$$
d_0=4(d^*)^2.
$$
We prove the two following properties.

 For any $d_0>0$
\begin{equation}\label{Eq_a}
\left\{
\Espace
\begin{array}{l}
\exists  \; p_0(d_0)>0, \; (p_i)_{i \geq 1}, \; T_0(d_0)>0  \; \textrm{ such that for any } l\leq k , \\
\P\left(l_0(k+1)=l \; | \; l_0(k)=l \right) \geq p_{k-l}, \textrm{ for any } T\geq T_0(d_0), \\
1-p_i \leq  e^{-i T}, \; i \geq 1,
\end{array}
\right.
\end{equation}
and, for any $(R_{0},d_0)$ sufficiently large,
\begin{equation}\label{Eq_b}
\left\{
\Espace
\begin{array}{l}
\exists \; T^*(R_0) >0  \textrm{ and } p_{-1}>0 \textrm{ such that for any }T\geq T^*(R_0)\\
\P\left(l_0(k+1)=k+1 \; | \; l_0(k)=\infty ,\; \mathcal H_k \leq R_0 \right)
 \geq p_{-1},
\end{array}
\right.
\end{equation}
where
$$
\mathcal H _k=\abs{u_1(kT)}^2+\abs{u_2(kT)}^2.
$$
\eqref{Eq_a} states that the probability that two solutions decouples at $kT$ is very small, \eqref{Eq_b} states that,
 inside a ball, the probability that two solutions get coupled at $(k+1)T$ is uniformly bounded below.

In the particular case where $\sig_l(x)$ does not depend on $x$ and where $K_0=0$, one can apply a similar proof
as in section 1.1 to establish a result closely related to  \eqref{Eq_a}, \eqref{Eq_b}. This technic has been developped in \cite{KS3}. But
it does not seem to work in the general case.

 Consequently, we use some tools developped in \cite{Matt} to establish  \eqref{Eq_a}, \eqref{Eq_b}.
 Note that in \eqref{Eq_b}, we use only starting points in a ball of radius $R_0$.  This is due to the fact that to prove
\eqref{Eq_b}, we need to estimate some terms which cannot be controlled on $\R^2$ but only inside
a ball. This further difficulty is due to the fact that contrary to the simple example of section 1.1, we work on an
unbounded phase space and is overcomed thanks to another ingredient which is the so-called Lyapunov structure. It allows
the control of the probability to enter the ball of radius $R_0$. In our example, it is an easy consequence of
\eqref{Eq_un_neuf_bis}ii). More precisely, we use the property that for any solution $u(\cdot,u_0)$
\begin{equation}\label{Eq_c}
\Espace
\left \{
\begin{array}{lcl}
\E\abs{u(t,u_0)}^2 &\leq & e^{-2t}\abs{u_0}^2+\frac{K_1}{2}, \\
\E\left(\abs{u(\tau',u_0)}^4 1_{\tau'<\infty}\right) & \leq &
K' \left( \abs{u_0}^4 + 1+\E\left(\tau' 1_{\tau'<\infty}\right) \right),
\end{array}
\right.
\end{equation}
 for any stopping times $\tau'$.

The following Proposition is a consequence of Theorem \ref{Th_Theo} given in a more general setting below.
\begin{Proposition}\label{Prop_theo}
If there exists a coupling of $\Dr(u(\cdot,u_0^i),W)$ such that
 \eqref{Eq_a}, \eqref{Eq_b} are satisfied, then \eqref{Eq_etoile} is true. Thus there
 exists a unique invariant measure $\nu$ of $(\mathcal P_t)_t$. Moreover there exist $C$ and $\alpha$ such that
$$
\norm{\mathcal P_t \mu -\nu}_{Lip_b(\R^2)}^*\leq Ce^{-\alpha t}\left( 1+\int_{\R^2}\abs{u}d\mu(u) \right).
$$
\end{Proposition}

To obtain \eqref{Eq_a} and \eqref{Eq_b}, we introduce three more ingredients. First in order
to build a coupling $((u_1,W_1),(u_2,W_2))$ such that
$((X_1,\eta_1),(X_2,\eta_2))$ is a maximal coupling, we use the following results contained in \cite{Matt}, although not
 explicitly stated. Its proof is postponed to the appendix.
\begin{Proposition}\label{Prop_Matt}
Let $E$ and $F$ be two polish spaces, $f_0:E\to F$ be a measurable map  and $(\mu_1,\mu_2)$ be two probability measures
on $E$. We set
$$
\nu_i=f_0^* \mu_i, \quad i=1,2.
$$
Then there exist a coupling $(V_1,V_2)$ of $(\mu_1,\mu_2)$ such that $(f_0(V_1),f_0(V_2))$ is a maximal coupling of
 $(\nu_1,\nu_2)$.
\end{Proposition}
We also remark that given $(X,\eta)$ on $[0,T]$, there exists a
 unique solution $Y(\cdot,u_0)$ of
$$
d Y+ 2 Y dt +g(X,Y)dt  = \sig_h(X) d\eta ,\quad Y(0,u_0)=y_0.
$$
We set
$$
Y(\cdot,u_0)=\Phi(X,\eta,u_0)(\cdot).
$$
It is easy to see that $Y$ is adapted to the filtration associated to $\eta$ and $X$.

Proposition \ref{Prop_Foais_Ex} implies that for any given $(X,\eta)$
\begin{equation}\label{Eq_dbis}
\abs{\Phi(X,\eta,u_0^1)(t)-\Phi(X,\eta,u_0^2)(t)}\leq e^{-t}\abs{u_0^1-u_0^2}.
\end{equation}

Then we rewrite the equation for $X$ as follows
\begin{equation}\label{Eq_d}
\Espace
\left\{
\Espace \begin{array}{lcl}
d X+ 2 X dt +f(X,\Phi(X,\eta,u_0))dt & = & \sig_l(X) d\beta , \\
\lefteqn{ X(0)=x_0.}
\end{array}
\right.
\end{equation}
The Girsanov formula can then be used on \eqref{Eq_d} as in section 1.1.

We finally remark that by induction, it suffices to construct a probability space $(\Omega_0,\F_0,\P_0)$ and two measurable couples of functions
 $(\omega_0,u_0^1,u_0^2)\to(V_i(\cdot,u_0^1,u_0^2))_{i=1,2}$ and  $(V_i'(\cdot,u_0^1,u_0^2))_{i=1,2}$ and  such that,
 for any
$(u_0^1,u_0^2)$, $(V_i(\cdot,u_0^1,u_0^2))_{i=1,2}$ and $(V_i'(\cdot,u_0^1,u_0^2))_{i=1,2}$ are two couplings of
$(\Dr(u(\cdot,u_0^i),W))_{i=1,2}$ on $[0,T]$. Indeed, we first set
$$
u_i(0)=u_0^i,\quad W_i(0)=0, \quad i=1,2.
$$
Assuming that we have built $(u_i,W_i)_{i=1,2}$ on $[0,kT]$, then we take $(V_i)_i$ and $(V'_i)_i$ as above independant
 of $(u_i,W_i)_{i=1,2}$ on $[0,kT]$ and set
\begin{equation}\label{Eq_m_bis}
\Espace
\left(u_i(kT+t),W_i(kT+t)\right) =\left \{
\Espace \begin{array}{ll}
V_i(t,u_1(kT),u_2(kT)) & \textrm{ if } l_0(k) \leq k,\\
V'_i(t,u_1(kT),u_2(kT))& \textrm{ if } l_0(k) =\infty,
\end{array}
\right.
\end{equation}
for any $t \in [0,T]$.

{\bf Proof of \eqref{Eq_a}.}

To build $(V_i(\cdot,u_0^1,u_0^2))_{i=1,2}$, we apply Proposition \ref{Prop_Matt} to $E=C((0,T);\R^2)^2$, $F=C((0,T);\R)^2$,
$$
f_0\left(u,W\right)= (X,\eta), \; \textrm{ where } u=\left( \begin{array}{c} X \\ Y \end {array} \right), \;
 W=\left( \begin{array}{c} \beta \\ \eta \end {array} \right),
$$
and to
$$
\mu_i=\Dr(u(\cdot,u_0^i),W), \quad \textrm{ on } [0,T].
$$
Remark that if we set $\nu_i=f_0^*\mu_i$, we obtain
$$
\nu_i=\Dr(X(\cdot,u_0^i),\eta), \quad \textrm{ on } [0,T].
$$
We write
$$
(Z_i,\xi_i)=f_0(V_i), \quad i=1,2.
$$

Then $(V_i(\cdot,u_0^1,u_0^2))_{i=1,2}$ is a coupling of $(\mu_1,\mu_2)$ such that
$((Z_i,\xi_i)(\cdot,u_0^1,u_0^2))_{i=1,2}$ is a maximal coupling of $(\nu_1,\nu_2)$.

 We first use a Girsanov formula to estimate $I_p$, where
$$
I_p=\int_F \left ( \frac{d\nu_2}{d \nu_1} \right )^{p+1} d\nu_2.
$$
Then, using Lemma \ref{lem_norm_var_sens}, we establish \eqref{Eq_a}.

We consider a couple $(u_i,W_i)_{i=1,2}$ consisting of
two solutions of \eqref{Eq_un_neuf} on $[0,kT]$. From now on, we are only concerned with a trajectory of $(u_i,W_i)_{i=1,2}$
 such that $l_0(k)=l\leq k$.  We set
$$
x=X_1(kT)=X_2(kT), \quad y_i=Y_i(kT), \quad i=1,2.
$$
Let $(\beta,\xi)$ be a two-dimensionnal brownian motion
 defined on a probability space $(\Omega,\F,\P)$. We denote by $Z$ the unique solution of
\begin{equation}\label{Eq_i_un}
\Espace
\left\{
\begin{array}{lcl}
d Z+ 2 Z dt +f(Z,\Phi(Z(\cdot),\xi(\cdot),(x,y_1)))dt & = & \sig_l(Z) d\beta , \\
\lefteqn{ Z(0)=x.}
\end{array}
\right.
\end{equation}
Taking into account \eqref{Eq_i_un}, we obtain that $\nu_1$ is the distribution of
$(Z,\xi)$ under the probability $\P$.

\noindent We set $\widetilde{\beta}(t)=\beta(t)+\int_0^t d(s) dt $ where
\begin{equation}\label{Eq_j}
 d(t)=\frac{1}{\sig_l(Z(t))}\left(f(Z(t),\Phi(Z,\xi,(x,y_2))(t))-f(Z(t),\Phi(Z,\xi,(x,y_1))(t))\right).
\end{equation}
Then $Z$ is a solution of
\begin{equation}\label{Eq_i_deux}
\Espace
\left\{
\begin{array}{lcl}
d Z+ 2 Z dt +f(Z,\Phi(Z(\cdot),\xi(\cdot),(x,y_2)))dt & = & \sig_l(Z) d\widetilde{\beta} , \\
\lefteqn{ Z(0)=x.}
\end{array}
\right.
\end{equation}
Since $f$ is bounded and $\sig_l$ is bounded below, then $d$ is uniformly bounded. Hence, the Novikov condition is satisfied and the Girsanov formula
can be applied. Then we set
$$
d\widetilde{\P}= \exp\left(\int_0^T d(s) dW(s)-\frac{1}{2} \int_0^T \abs{d(s)}^2 dt  \right)d\P
$$
 We deduce from the Girsanov formula that $\widetilde{\P}$ is a probability under which $(\widetilde{\beta},\xi)$ is a
 brownian motion and since $Z$ is a solution of \eqref{Eq_i_deux}, then the law of $(Z,\xi)$ under $\widetilde{\P}$ is $\nu_2$.
 Moreover
\begin{equation}\label{Eq_jbis}
I_p\leq \E
\exp\left( c_p \int_0^T \abs{d(s)}^2 dt \right).
\end{equation}
Since $f$ is Lipschitz, then we infer from \eqref{Eq_j} and \eqref{Eq_un_neuf_ter}i) that
$$
\abs{d(t)} \leq \sig_0^{-1}\abs{f}_{Lip} \abs{\Phi(Z(\cdot),\xi(\cdot),(x,y_1))(t)-\Phi(Z(\cdot),\xi(\cdot),(x,y_2))(t)}.
$$
Now we use the Foias-Prodi estimate. Applying \eqref{Eq_abis} and \eqref{Eq_dbis}, it follows from $l_0(k)=l$ that
$$
\abs{d(t)}^2 \leq d_0 \sig_0^{-2}\abs{f}_{Lip}^2  \exp\left( -2(k-l)T \right).
$$
Then it follows that
\begin{equation}\label{Eq_k}
I_p\leq
\exp\left( c_p \sig_0^{-2}d_0 \abs{f}_{Lip}^2  e^{ -2(k-l)T } \right).
\end{equation}
Note that
$$
\norm{\nu_1-\nu_2}_{var}=\int_F \abs{ \frac{d \nu_2}{d \nu_1}  -1}d\nu_2
\leq \sqrt{\int \left ( \frac{d \nu_2}{d \nu_1} \right )^{2} d\nu_2-1}.
$$
We infer from \eqref{Eq_k} that, for $T\geq T_0(d_0)=( \sig_0^{-2}c_p d_0 \abs{f}_{Lip}^2)^{-1}$,
$$
\norm{\nu_1-\nu_2}_{var}\leq e^{ -(k-l)T}.
$$
 Applying Lemma \ref{lem_norm_var_sens} to the maximal coupling $(Z_1,Z_2)_{i=1,2}$ of $(\nu_1,\nu_2)$ gives
\begin{equation}\label{Eq_l}
\P\left( (Z_1,\xi_1) \not = (Z_2,\xi_2) \right) \leq \norm{\nu_1-\nu_2}_{var}\leq e^{ -(k-l)T}.
\end{equation}
Using \eqref{Eq_m_bis} and \eqref{Eq_l}, we obtain that  on $l_0(k)=l$
$$
\P\left( (X_1,\eta_1) \not = (X_2,\eta_2) \textrm{ on } [kT,(k+1)T]  \,\left| \, \F_{kT} \right.\right)
\leq  e^{ -(k-l)T}.
$$
Noticing that
$$
\{l_0(k+1)=l\}=\{l_0(k)=l\}\cap\{(X_1,\eta_1)  = (X_2,\eta_2) \textrm{ on } [kT,(k+1)T]\}.
$$
and integrating over $l_0(k)=l$ gives for $T\geq T_0(d)$ and for $k>l$
\begin{equation}\label{Eq_e}
\P\left( l_0(k+1)\not=l \,|\,l_0(k)=l\right)
\leq e^{ -(k-l)T}.
\end{equation}

Now, it remains to consider the case $k=l$, we apply Lemmas  \ref{lem_norm_var_sens} and
 \ref{lem_tech_coup_inf} to
 $(Z_i,\xi_i)_{i=1,2}$ which gives
$$
\P\left( (Z_1,\xi_1)=(Z_2,\xi_2) \right)= \left( \nu_1\wedge\nu_2 \right)(F)\geq \left(1-\frac{1}{p}\right) (pI_p)^{-\frac{1}{p-1}}.
$$
Applying \eqref{Eq_jbis} and fixing $p>1$, we obtain
\begin{equation}\label{Eq_ebis}
\P\left( (Z_1,\xi_1)=(Z_2,\xi_2) \right)\geq p_0(d_0)= \left(1-\frac{1}{p}\right)p^{-\frac{1}{p-1}}
\exp\left( - c_p d_0 \abs{f}_{Lip}^2  \right).
\end{equation}
To conclude, we notice that \eqref{Eq_e} and \eqref{Eq_ebis} imply \eqref{Eq_a}.

{\bf Proof of \eqref{Eq_b}.}

Assume that we have $d_0>0$, $\widetilde p>0$, $T_1>0$, $R_1>4 K_1$ and  a coupling $(\widetilde V_i(\cdot,u_0^1,u_0^2))_{i=1,2}$ of
$(  \mu_1,  \mu_2)$, where
$$
  \mu_i=\Dr (u(\cdot,u_0^i),W), \; \textrm{ on } [0,T_1], \quad i= 1,2,
$$
and such that for any $(u_0^1,u_0^2)$ which satisfies $\abs{u_0^1}^2+\abs{u_0^2}^2\leq R_1$
\begin{equation}\label{Eq_b_b}
\P\left(  Z_1(T_1,u_0^1,u_0^2)=  Z_2(T_1,u_0^1,u_0^2),\;\sum_{i=1}^2\abs{  u_i(T_1,u_0^1,u_0^2)}^2\leq d_0
 \right)\geq \widetilde p,
\end{equation}
where
$$
\widetilde V_i(\cdot,u_0^1,u_0^2)=\left(   u_i(\cdot,u_0^1,u_0^2),   W_i(\cdot,u_0^1,u_0^2) \right),
\;   u_i(\cdot,u_0^1,u_0^2)=\left(\begin{array}{c}   Z_i\\  G_i \end{array}\right), \; i=1,2.
$$
By applying the Lyapunov structure \eqref{Eq_c}, we obtain that for any $\theta\geq T_2(R_0,R_1)$
\begin{equation}\label{Eq_b_c}
\P\left(\abs{u(\theta,u_0)}^2\geq \frac{R_1}{2}\right)\leq \frac{1}{4}, \; \textrm{ for any } u_0 \textrm{ such that }
 \abs{u_0}^2 \leq \frac{R_0}{2}.
\end{equation}
In order to build $(V'_1,V'_2)$ such that \eqref{Eq_b} happens, we set $T^*(R_0)=T_1+T_2(R_0)$ and for any $T\geq T^*(R_0)$,
 we set $\theta=T-T_1$ and we remark that $\theta\geq T_2(R_0)$. Then we construct the trivial coupling $(V"_1,V"_2)$ on $[0,\theta]$. Finally, we consider
 $(\widetilde V_1,\widetilde V_2)$ as above independant of $(V"_1,V"_2)$ and we set
$$
\Espace
V'_i(t,u_0^1,u_0^2)=\left\{
\begin{array}{ll}
V"_i(t,u_0^1,u_0^2) & \textrm{ if } t\leq \theta, \\
\widetilde V_i(t-\theta,V"_1(\theta,u_0^1,u_0^2),V"_2(\theta,u_0^1,u_0^2)) & \textrm{ if } t\geq \theta.
\end{array}
\right.
$$
Combining \eqref{Eq_b_b} and \eqref{Eq_b_c}, we obtain \eqref{Eq_b}
  with $p_{-1}=\frac{1}{2}\widetilde p$.

To build $(\widetilde V_i(\cdot,u_0^1,u_0^2))_{i=1,2}$, we apply Proposition \ref{Prop_Matt} to $E=C((0,T_1);\R^2)^2$,
$F=\R$,
$$
f_0\left(u,W\right)= X(T_1), \; \textrm{ where } u=\left( \begin{array}{c} X \\ Y \end {array} \right), \;
 W=\left( \begin{array}{c} \beta \\ \eta \end {array} \right),
$$
and to $(  \mu_1, \mu_2)$.
Remark that if we set $ \nu_i=f_0^* \mu_i$, we obtain
$$
 \nu_i=\Dr(X(T_1,u_0^i)).
$$

Then $(\widetilde V_i(\cdot,u_0^1,u_0^2))_{i=1,2}$ is a coupling of $( \mu_1, \mu_2)$ such that
$(Z_i(T_1,u_0^1,u_0^2))_{i=1,2}$ is a maximal coupling of $( \nu_1, \nu_2)$.

 Now we notice that if we have
$(\hat{\nu}_1,\hat \nu_2)$ two equivalent measures such that $ \nu_i$ is equivalent to $\hat{\nu}_i$ for $i=1,2$, then by applying
 two Schwartz inequality, we obtain that
\begin{equation}\label{Eq_b_a}
I_p \leq \left( J_{2p+2}^1\right)^{\frac{1}{2}}
\left(J_{4p}^2\right)^{\frac{1}{4}}\left(\hat I_{4p+2}\right)^{\frac{1}{4}},
\end{equation}
where $A=[-d_1,d_1]$ and
$$
\Espace
\begin{array}{rclrcl}
I_p &=& \int_A \left ( \frac{d  \nu_1}{d  \nu_2} \right )^{p+1} d \nu_2,&
J_p^1 &=& \int_A \left ( \frac{d  \nu_1}{d \hat\nu_1} \right )^{p} d\hat\nu_1,\\
\hat I_p &=& \int_A \left ( \frac{d \hat\nu_1}{d \hat\nu_2} \right )^{p} d\hat\nu_2,&
J_p^2 &=& \int_A \left ( \frac{d \hat\nu_2}{d  \nu_2} \right )^{p} d\hat\nu_2
\end{array}
$$
Recall that $Z_i$ the unique solution of
\begin{equation}\label{Eq_b_d_un}
\Espace
\left\{
\begin{array}{lcl}
d Z_i+ 2 Z_i dt +f(Z_i,\Phi(Z_i(\cdot),\xi_i(\cdot),u_0^i))dt & = & \sig_l(Z_i) d\beta_i , \\
\lefteqn{ Z_i(0)=x_0^i.}
\end{array}
\right.
\end{equation}

\noindent We set $\widetilde{\beta}_i(t)=\beta_i(t)+\int_0^t d_i(s) dt $ where
\begin{equation}\label{Eq_b_d_j}
 d_i(t)=-\frac{1}{\sig_l(Z_i(t))}f(Z_i(t),\Phi(Z_i(\cdot),\xi_i(\cdot),u_0^i)(t)).
\end{equation}
Then $Z_i$ is a solution of
\begin{equation}\label{Eq_b_d_i_deux}
\Espace
\left\{
\begin{array}{lcl}
d Z_i+ 2 Z_i dt  & = & \sig_l(Z_i) d\widetilde{\beta}_i , \\
\lefteqn{ Z_i(0)=x_0^i.}
\end{array}
\right.
\end{equation}
Since $f$ is bounded and $\sig_l$ is bounded below, then $d_i$ is uniformly bounded. Hence, the Novikov condition is satisfied and the Girsanov formula
can be applied. Then we set
$$
d\widetilde{\P}_i= \exp\left(\int_0^T d_i(s) dW(s)-\frac{1}{2} \int_0^T \abs{d_i(s)}^2 dt  \right)d\P
$$
 We deduce from the Girsanov formula that $\widetilde{\P}_i$ is a probability under which $(\widetilde{\beta}_i,\xi_i)$ is a
 brownian motion. We denote by $\hat \nu_i$ the law of $Z_i(T_1)$ under $\widetilde{\P}_i$.
 Moreover
\begin{equation}\label{Eq_b_d_jbis}
J_p^i\leq
\exp\left( c_p \int_0^T \abs{d_i(s)}^2 dt \right)\leq \exp\left( c_p \sig_0^{-2} \abs{f}_\infty^2 \right).
\end{equation}
It is classical that since $\sig_l$ is bounded below, then $\hat\nu_i$ has a density $q(x_0^i,z)$ with respect to lebesgue measure $dz$, that $q$ is continuous with
 respect to the couple
$(x_0^i,z)$, where $x_0^i$ is the initial value and where $z$ is the target value and that $q>0$. Then, we can
 bound $q$ and $q^{-1}$ uniformly on
$\abs{x_0^i}\leq R_1$ and $z\in A=[-d_1,d_1]$, which allows us to bound $\hat I_p$ and then $I_p$. Actually:
\begin{equation}\label{Eq_b_e}
I_p\leq C'(p,d_1,T_1,R_1)<\infty.
\end{equation}
Now we apply Lemmas  \ref{lem_tech_coup_inf} and  \ref{lem_norm_var_sens}:
\begin{equation}\label{Eq_b_ebis}
\P\left(  Z_1(T_1)=  Z_2(T_1),\; \abs{  Z_1(T_1)}\leq d_1\right)\geq
\left(1-\frac{1}{p}\right)p^{-\frac{1}{p-1}}I_p^{-\frac{1}{p-1}}\nu_1([-d_1,d_1])^{\frac{p}{p-1}}.
\end{equation}
If we fix $d_1>4K_1$, then we obtain from the Lyapunov structure \eqref{Eq_c} that there exists $T_1=T_1(R_1,d_1)$ such
that
\begin{equation}\label{Eq_b_eter}
\nu_1([-d_1,d_1])\geq \frac{1}{2}.
\end{equation}
Combining \eqref{Eq_b_e}, \eqref{Eq_b_ebis} and \eqref{Eq_b_eter} gives
\begin{equation}\label{Eq_b_f}
\P\left(  Z_1(T_1)=  Z_2(T_1),\; \abs{  Z_1(T_1)}\leq d_1\right)\geq
 C(p,d_1,T_1,R_1)>0.
\end{equation}
Note that
\begin{equation}\label{Eq_b_g}
\Espace
\begin{array}{ll}
\P(&  Z_1(T_1)=\;  Z_2(T_1),\; \abs{  u_i(T_1)}\leq d_1+d_2,\; i=1,2)\geq\\
&\P\left(  Z_1(T_1)=  Z_2(T_1),\; \abs{  Z_1(T_1)}\leq d_1\right)
-\sum_{i=1}^2\P\left(\abs{  u_i(T_1)}\geq d_2\right).
\end{array}
\end{equation}
Using the Lyapunov structure \eqref{Eq_c}, we obtain that
\begin{equation}\label{Eq_b_h}
\P\left(\abs{  u_i(T_1)}\geq d_2\right)\leq \frac{R_1+K_1}{d_2^2}.
\end{equation}
Combining \eqref{Eq_b_f}, \eqref{Eq_b_g} and \eqref{Eq_b_h}, we can choose $d_2$ sufficiently high such that, by setting
$d^*=d_1+d_2$, $d_0=(2d^*)^2$ and $\hat p=\frac{1}{2}C(2,d_1,T_1,R_1)$, \eqref{Eq_b_b} holds.
\subsection{Abstract Result}
\

We now state and prove an abstract result which allows to reduce the proof of exponential convergence to equilibrium
 to the verification of some conditions, as was done in the previous section.

This result is closely related to the abstract result of \cite{Matt}. Our proof has some similarity with the one in
the reference but, in fact, is closer to arguments used in \cite{S}. Our abstract result could be used in articles
\cite{K}, \cite{KS}, \cite{KS2}, \cite{KS3} and \cite{Matt} to conclude.

In fact, in \cite{Matt} a family $(r_k,s_k)$ of subprobability are used, whereas in \cite{K}, \cite{KS3} a family of
subsets $Q(l,k)$ are introduced. Here, we use a random integer valued process $l_0(k)$. The three
points of view are equivalent, the correspondance is given by
$$
s_{k+1}=\P\left(\{l_0(k+1)=1\}\cap \cdot\right),\quad r_{k+1}=\P\left(\{l_0(k+1)=1\}^c\cap \cdot\right),
$$
and
$$
Q(l,k)=\{l_0(k)=l\}.
$$
The result has already been applied in section 1.2, the function used below is
$$
\Hcal(u_0)=\abs{u_0}^2,
$$
in this example. In fact, in most of the application and in particular for the CGL equation in the first case
 treated below, $\Hcal$ wil be the square of the norm.
We are concerned with
$v(\cdot,(u_0,W_0))=(u(\cdot,u_0),W(\cdot,W_0))$, a couple of strongly Markovian process defined on polish spaces $(E,d_E)$ and
$(F,d_F)$.
 We denote by $(\Pcal_t)_{t\in I}$ the markovian transition semigroup of $u$, where $I=\R^+$ or $T\N=\{kT, \; k\in \N\}$.

We consider for any initial conditions $(v_0^1,v_0^2)$  a coupling $(v_1,v_2)$ of
$(\Dr(v(\cdot,v_0^1)),\Dr(v(\cdot,v_0^2)))$ and a random integer valued process $l_0:\N\to \N\cup \{\infty\}$ which has
 the following properties
\begin{equation}\label{Abstract_ater}
\Espace
\left\{
\begin{array}{l}
l_0(k+1)=l \textrm{ implies } l_0(k)=l, \textrm{ for any } l\leq k,\\
l_0(k)\in \{0,1,2, ...,k\}\cup\{\infty\},\\
l_0(k) \textrm{ depends only of } v_1|_{[0,kT]} \textrm{ and } v_2|_{[0,kT]},\\
l_0(k)=k\; \textrm{ implies } \Hcal_k \leq d_0,
\end{array}
\right.
\end{equation}
where
$$
\mathcal H_k=\mathcal H(u_1(kT))+\mathcal H(u_2(kT)), \quad \mathcal H:\, E \to \R^+.
$$
We write $v_i=(u_i,W_i)$.
 From now on we say that $(v_1,v_2)$ are coupled at $kT$ if $l_0(k)\leq k$, in other words
 if $l_0(k)\not = \infty$.

Now we see four conditions on the coupling. The first condition states that when $(v_1,v_2)$
 have been coupled for a long time then the probability that $(u_1,u_2)$ are close is high.
\begin{equation}\label{Abstract_abis}
\Espace
\left\{
\begin{array}{l}
\textrm{There exist $c_0$ and $\alpha_0>0$ such that }\\
\P\left(    d_E(u_1(t),u_2(t))> c_0 e^{-\alpha_0(t-lT)} \; \textrm{ and } \; l_0(k)=l\right)
 \leq c_0 e^{-\alpha_0(t-lT)},
\end{array}
\right.
\end{equation}
 for any  $t \in [lT,kT]\cap I$.

The following property states that the probability that two solutions decouples at $kT$ is very small
\begin{equation}\label{Abstract_a}
\left\{
\Espace
\begin{array}{l}
\textrm{There exist }  \;  (p_k)_{k \in \N},\; c_1>0,\; \alpha_1>0  \; \textrm{ such that} , \\
\P\left(l_0(k+1)=l \; | \; l_0(k)=l \right) \geq p_{k-l}, \textrm{ for any } l\leq k , \\
1-p_k \leq  c_1 e^{-\alpha_1 k T},\; p_k>0  \;\textrm{ for any }  k \in \N.
\end{array}
\right.
\end{equation}
Next condition states that,
 inside a ball, the probability that two solutions get coupled at $(k+1)T$ is uniformly bounded below.
\begin{equation}\label{Abstract_b}
\left\{
\Espace
\begin{array}{l}
\textrm{There exist }  p_{-1}>0, \; R_0>0 \textrm{ such that}\\
\P\left(l_0(k+1)=k+1 \; | \; l_0(k)=\infty ,\; \mathcal H_k \leq R_0 \right)
 \geq p_{-1}.
\end{array}
\right.
\end{equation}

The last ingredient is the so-called Lyapunov structure. It allows
the control of the probability to enter the ball of radius $R_0$. It states that there exists $\gamma> 1$, such that
 for any solution $v_0$
\begin{equation}\label{Abstract_c}
\Espace
\left \{
\begin{array}{lcl}
\E\mathcal H(v(t,v_0)) &\leq & e^{-\alpha_3t}\mathcal H(v_0)+\frac{K_1}{2}, \\
\E\left(\mathcal H(v(\tau',v_0))^{\gamma} 1_{\tau'<\infty}\right) & \leq &
K' \left(\mathcal H(v_0) + 1+\E\left(\tau' 1_{\tau'<\infty}\right) \right)^{\gamma},\\
\lefteqn{\textrm{for any stopping times $\tau'$ taking value in $\{kT,\; k\in \N\}\cup\{\infty\}$.}}
\end{array}
\right.
\end{equation}

The process $V=(v_1,v_2)$ is said to be $l_0$--Markovian if the laws of $V(kT+\cdot)$ and of $l_0(k+\cdot)-k$ on
 $\{l_0(k)\in \{k,\infty\}\}$
 conditionned by $\F_{kT}$ only depend on $V(kT)$ and are equal to the laws of $V(\cdot,V(kT))$ and $l_0$, respectively.

Notice that in the example of the previous section or in the CGL case below, the process $(u_i,W_i)_{i=1,2}$ is $l_0$--Markovian
 but not Markovian. However, in both cases, if we choose $d_0=R_0$, we can modify the coupling such that the couple
 is Markovian at discrete times $T\N=\{kT,\; k\in \N\}$. But it does not seem to be possible to modify the coupling
     to become Markovian at any times.

\begin{Theorem}\label{Th_Theo}
Assume that \eqref{Abstract_ater}, \eqref{Abstract_abis},
 \eqref{Abstract_a}, \eqref{Abstract_b} and
\eqref{Abstract_c} hold whith $R_0>4 K_1$ and $R_0\geq d_0$ and that $V=(v_1,v_2)$ is $l_0$--Markovian. Then
 there exist $\alpha_4>0$ and $c_4>0$  such that
\begin{equation}\label{Abs_a}
\P\left( d_E(u_1(t),u_2(t))>c_3 e^{-\alpha_4 t} \right)\leq c_3 e^{-\alpha_4 t}\left( 1+\Hcal(u_0^1)+\Hcal(u_0^2) \right).
\end{equation}
Moreover there exists a unique stationnary probability mesure $\nu$ of $(\Pcal_t)_{t\in I}$ on $E$.
It satisfies,
\begin{equation}\label{Eq1.49bis}
\int_E \mathcal H(u)d\nu(u) \leq \frac{K_1}{2},
\end{equation}
and there exists $c_4>0$  such that for any $\mu \in \Pcal(E)$
\begin{equation}\label{Un_49_bis}
\abs{\Pcal^*_t\mu-\nu}_{Lip_b(E)}^*\leq c_4 e^{-\alpha_4 t}\left(1+ \int_E \mathcal H(u)d\mu(u) \right).
\end{equation}
\end{Theorem}

Proposition \ref{Prop_theo} is an easy consequence of Theorem \ref{Th_Theo}. Actually \eqref{Abstract_ater} is clear and
 \eqref{Abstract_abis} and \eqref{Abstract_c} are consequence of \eqref{Eq_abis} and \eqref{Eq_c} if
$R_0\geq d_0$. Finally,
since, for any $(R_0,d_0,T)$ sufficiently high, there exists a coupling such that \eqref{Eq_a} and \eqref{Eq_b} hold,
we can choose $(R_0,d_0,T)$ such that all our assumptions are true.

\begin{Remark}
 Inequality \eqref{Un_49_bis} means that for any $f \in Lip_b(E)$ and any $u_0 \in E$
$$
\abs{\E f(u(t,u_0))-\int_E f(u)d\nu (u)} \leq c_4 \abs{f}_{lip_b(E)} e^{-\alpha_4 t}(1+\Hcal(u_0) ).
$$
\end{Remark}

\subsection{Proof of Theorem \ref{Th_Theo}}
\

{\bf Reformulation of the problem}

We rewrite our problem in the form on a exponential estimate.

As in the example, it is sufficient to establish \eqref{Abs_a}. Then \eqref{Eq1.49bis} is a simple consequence of
\eqref{Abstract_c} and \eqref{Un_49_bis} follows from \eqref{Eq1_15bis}.
Assume that $t> 8T$. We denote by $k$ the unique integer such that \mbox{$t\in (2(k-1)T,2kT]$}. Notice that
$$
\Espace
\begin{array}{l}
\lefteqn{\P( d_E(u_1(t),u_2(t))> c_0 e^{-\alpha_0 (t-(k-1)T)})}\\
\quad \quad \quad \leq \P\left(l_0(2k)\geq k\right)+
 \P\left( d_E(u_1(t),u_2(t))>c_0 e^{-\alpha_0 (t-(k-1)T)}  \; \textrm{ and } \; l_0(2k)< k \right).
\end{array}
$$
Thus applying \eqref{Abstract_abis}, using $2(t-(k-1)T)>t$, it follows
\begin{equation}\label{Abs_b}
\Espace
\begin{array}{l}
\P\left( d_E(u_1(t),u_2(t))>c_0 \exp\left({-\frac{\alpha_0}{2} t}\right) \right)\leq
\P\left(l_0(2k)\geq k\right)+
c_0 \exp\left({-\frac{\alpha_0}{2} t}\right).
\end{array}
\end{equation}
In order to estimate $\P\left(l_0(2k)\geq k\right)$, we introduce the following notation
$$
l_0(\infty)=\lim \sup l_0.
$$
Taking into account \eqref{Abstract_ater}, we obtain that for $l<\infty$
$$\{l_0(\infty)=l\}=\{ l_0(k)=l, \textrm{ for any } k\geq l\}.$$
We deduce
\begin{equation}\label{Abs_c}
\P\left(l_0(2k)\geq k\right)\leq \P\left(l_0(\infty)\geq k\right).
\end{equation}
Taking into account \eqref{Abs_b},  \eqref{Abs_c}  and using a Chebyshev inequality,
it is sufficient to obtain that there exist $c_5>0$ and $\delta>0$ such that
\begin{equation}\label{Abs_d}
\E \left( \exp\left(\delta l_0(\infty)\right)\right)\leq  c_5 \left(1+\Hcal(u_0^1)+\Hcal(u_0^2) \right).
\end{equation}
Then \eqref{Abs_a} follows with
$$
\alpha_4=\min\left\{\frac{\alpha_0}{2},\frac{\delta}{2T}\right\}.
$$

{\bf Definition of a sequence of stopping times}

Using the Lyapunov structure \eqref{Abstract_c}, we prove at the end this subsection that there exist $\delta_0>0$
and $c_6>0$ such that
\begin{equation}\label{Abs_e}
\E \left(\exp\left(\delta_0 \tau \right)\right)\leq  c_6 \left(1+\Hcal(u_0^1)+\Hcal(u_0^2) \right),
\end{equation}
where
$$
\tau=\min\left\{t\in T\N \,|\, \Hcal(u_1(t))+\Hcal(u_2(t))\leq R_0 \right\}.
$$

We set
$$
\hat \sig =\min \left\{ k\in \N^* \,|\, l_0(k)>1  \right\},\quad
\sig      =\hat \sig T.
$$
Clearly $\hat \sig = 1$ if the two solutions do not get coupled at time $0$ or $T$. Otherwise, they get coupled at $0$ or
$T$ and remain coupled until $\sig$.

Let us assume for the moment that if $\Hcal_0\leq R_0$, then
\begin{equation}\label{Abs_f}
\Espace
\left\{
\begin{array}{l}
\E \left (\exp\left(\delta_1 \sig \right)1_{\sig<\infty}\right)\leq  c_7 ,\\
\P\left( \sig=\infty \right)\geq p_\infty >0.
\end{array}
\right.
\end{equation}
The proof is given after the proof of \eqref{Abs_e} at the end of this subsection.

Now we build a sequence of stopping times
$$
\Espace
\begin{array}{lcllcl}
\tau_0&=&\tau,\\
\hat \sig_{k+1} &=&\min \left\{ l\in \N^* \,|\, lT>\tau_k \textrm { and }l_0(l) T>\tau_k+T  \right\},&
 \sig_{k+1}&=&\hat \sig_{k+1} \times T\\
\tau_{k+1}       &= & \sig_{k+1}+\tau o \theta_{\sig_{k+1}},
\end{array}
$$
 where $(\theta_t)_t$ is the shift operator.
The idea is the following. We wait the time $\tau_k$ to enter the ball of radius $R_0$. Then, if we do not start
 coupling at time $\tau_k$, we try to couple at time $\tau_k+T$. If we fail to start coupling at time $\tau_k$ or
$\tau_k+T$ we set $\sig_k=\tau_k+T$ else we set $\sig_k$ the time the coupling fails ($\sig_k=\infty$ if the
 coupling never fails). Then if $\sig_k<\infty$, we retry to enter the ball of radius $R_0$. The fact that $R_0\geq d_0$
 implies that  $l_0(\tau_k)\in \{\tau_k,\infty\}$.

The idea of the $l_0$--Markovian property is the following. Since $l_0(\tau_k)\in \{\tau_k,\infty\}$ and
$l_0(\sig_k)\in\{\sig_k,\infty\}$, when these stopping times are finite and  since these stopping times are taking value
 in  $T\N\cup\{\infty\}$, then the $l_0$--Markovian property implies the strong Markovian property when conditionning
 with respect to $\F_{\tau_k}$ or $\F_{\sig_k}$. Moreover, we infer from the $l_0$--Markovian property of $V$ that
$$
\sig_{k+1}=\tau_k+\sig o \theta_{\tau_k},
$$
which implies
$$
\tau_{k+1}=\tau_k+\rho o \theta_{\tau_k}, \; \textrm{ where } \rho=\sig + \tau o \theta_{\sig} .
$$

{\bf Exponential estimate on $\rho$}

Before concluding, we establish that  there exist $K$ such that for any $V_0$ such that $\Hcal_0\leq R_0$ and for any
$\delta_2\leq\frac{1}{\gamma'}\left( \delta_0\wedge\delta_1 \right)$
\begin{equation}\label{Abs_i}
\E_{V_0} \left(e^{\delta_2  \rho }1_{\rho<\infty}\right) \leq  K.
\end{equation}

Notice that for any $V_0$ such that $\Hcal_0\leq R_0$,
$$
\E_{V_0} \left(e^{\delta_2  \rho }1_{\rho<\infty}\right) =
\E_{V_0} \left(e^{\delta_2  \sig }1_{\sig<\infty}
\E \left(e^{\delta_2  \tau o \theta_{\sig} }1_{\tau o \theta_{\sig}<\infty} | \F_{\sig}\right)
\right) .
$$
Applying the $l_0$--Markovian property and \eqref{Abs_e}, we obtain
$$
\E \left(e^{\delta_2  \tau o \theta_{\sig} }1_{\tau o
\theta_{\sig}<\infty} | \F_{\sig}\right) \leq  c_6
\left(1+\Hcal(u_1(\sig))+\Hcal(u_2(\sig)) \right)1_{\sig<\infty},
$$
which implies
$$
\E_{V_0} \left(e^{\delta_2  \rho }1_{\rho<\infty}\right) \leq
c_6\E_{V_0} \left(e^{\delta_2  \sig }1_{\sig<\infty}
\left(1+\Hcal(u_1(\sig))+\Hcal(u_2(\sig)) \right)
\right) .
$$
An H\"older inequality gives
$$
\E_{V_0} \left(e^{\delta_2  \rho }1_{\rho<\infty}\right) \leq
c_6\left(\E_{V_0} e^{\gamma' \delta_2  \sig
}1_{\sig<\infty}\right)^{\frac{1}{\gamma'}}
\left(\E_{V_0}\left(1+\Hcal(u_1(\sig))+\Hcal(u_2(\sig))
\right)^\gamma1_{\sig<\infty}\right)^{\frac{1}{\gamma}}.
$$
Applying the Lyapunov structure \eqref{Abstract_c} and \eqref{Abs_f}, we obtain \eqref{Abs_i}.

{\bf Conclusion}

We remark that
$$
\E \left( e^{\delta_2  \tau_{k+1} }1_{\tau_{k+1}<\infty}\right) =
\E\left( e^{\delta_2  \tau_{k} }1_{\tau_{k}<\infty}
\E \left(e^{\delta_2  \rho o \theta_{\tau_k} }1_{\rho o \theta_{\tau_k}<\infty} | \F_{\tau_k}\right)\right) .
$$
Applying again the $l_0$-Markov property of $V$
\begin{equation}\label{Abs_g}
\E \left(e^{\delta_2  \tau_{k+1} }1_{\tau_{k+1}<\infty}\right) =
\E\left( e^{\delta_2  \tau_{k} }1_{\tau_{k}<\infty}
\E_{V(\tau_k    )} \left(e^{\delta_2  \rho }1_{\rho<\infty}\right)\right) .
\end{equation}
Iterating \eqref{Abs_g} by using \eqref{Abs_i} and \eqref{Abs_e}, we obtain
\begin{equation}\label{Abs_j}
\E e^{\delta_2  \tau_{n} }1_{\tau_{n}<\infty} \leq
 c_6 K^n\left(1+\Hcal(u_0^1)+\Hcal(u_0^2) \right).
\end{equation}
Using the second inequality of \eqref{Abs_f} and that $\tau<\infty$, we obtain from the $l_0$--Markov property that
\begin{equation}\label{Abs_k}
\P\left( k_0>n \right)\leq \left(1-p_\infty\right)^n,
\end{equation}
where
$$k_0=\inf\{k\in \N\,|\, \sig_{k+1}=\infty\}.$$
Then we obtain that $k_0<\infty$ almost surely and that
$$
l_0(\infty)\in\{\tau_{k_0},\tau_{k_0}+1\}.
$$
Therefore $l_0(\infty)<\infty$ almost surely and
$$
\E \exp\left( \frac{\delta_2}{p} l_0(\infty) \right) \leq\sum_{n=1}^\infty \E e^{\frac{\delta_2}{p}  (\tau_{n}+1) }1_{k_0=n},
$$
which implies, by applying a H\"older inequality,
$$
\E \exp\left( \frac{\delta_2}{p} l_0(\infty) \right) \leq e^{\frac{\delta_2}{p}  }
\sum_{n=1}^\infty \left(\E e^{\delta_2  \tau_{n} }1_{\tau_n\leq \infty}\right)^{\frac{1}{p}}
\left( \P\left( k_0=n \right)  \right)^{\frac{1}{p'}}.
$$
Applying \eqref{Abs_j} and \eqref{Abs_k}, we obtain
$$
\E \exp\left( \frac{\delta_2}{p} l_0(\infty) \right) \leq
c_6e^{\frac{\delta_2}{p}  }\left(\sum_{n=1}^\infty \left(K^{\frac{1}{p}} (1-p_\infty)^{\frac{1}{p'}}\right)^n\right)
\left(1+\Hcal(u_0^1)+\Hcal(u_0^2) \right)^\frac{1}{p}.
$$
Choosing $p$ such that $K^{\frac{1}{p}} (1-p_{\infty})^{\frac{1}{p'}}<1$ and setting $\delta= \frac{\delta_2}{p}$, we obtain \eqref{Abs_d}

{\bf Proof of \eqref{Abs_e}}

Let $N$ be an integer such that
$$
e^{-\alpha_3 N T}\leq \frac{1}{8}.
$$
We fix $i\in \{1,2\}$ and  set
$$
B_k=\left\{ \Hcal\left(u_i\left(  j N T  \right)\right)\geq 2 K_1, \textrm{ for any } j\leq k \right\},
\quad C_k=\left\{ \Hcal\left(u_i\left(  k N T  \right)\right)\geq 2 K_1\right\}.
$$
Combining the Markov property of $u_i$ and  the Lyapunov structure \eqref{Abstract_c}, we obtain
\begin{equation}\label{Abs_e_a}
\E \left(\Hcal(u_i((k+1)NT))|\F_{kNT}\right)\leq  \frac{1}{4}\Hcal(u_i(kNT))+\frac{K_1}{2} .
\end{equation}
Hence, applying a Chebyshev inequality, it follows that
\begin{equation}\label{Abs_e_b}
\P \left(C_{k+1}|\F_{kNT}\right)\leq  \frac{1}{8 K_1}\Hcal(u_i(kNT))+\frac{1}{4} .
\end{equation}
Integrating \eqref{Abs_e_a}, \eqref{Abs_e_b} over $B_k$, we obtain that
\begin{equation}\label{Abs_e_c}
\Espace
\left(
\begin{array}{c}
\E \left( \Hcal(u_i((k+1)NT))1_{B_{k+1}} \right)\\
\P\left( B_{k+1} \right)
\end{array}
\right)
\leq  A
\left(
\begin{array}{c}
\E \left( \Hcal(u_i(kNT))1_{B_{k}} \right)\\
\P\left( B_{k} \right)
\end{array}
\right),
\end{equation}
where
$$
\Espace
A=
\left(
\begin{array}{cc}
\frac{1}{4}     &    \frac{K_1}{2}   \\
\frac{1}{8K_1}  &    \frac{1}{4}
\end{array}
\right).
$$
Since the eigenvalues of $A$ are $0$ and $\frac{1}{2}$, we obtain that
$$
\P\left( B_{k} \right)\leq \frac{2}{K_1}\left(\frac{1}{2}\right)^k \left(1+\Hcal(u_0^i)\right).
$$
It follows from $R_0\geq 4 K_1$ that
$$
\P\left( \tau >k T \right)\leq c\exp\left(-\frac{k}{N}\ln 2\right) \left(1+\Hcal(u_0^i)\right).
$$
Hence, taking $\delta_0<\frac{\alpha_3}{3}$, we have established \eqref{Abs_e} .

{\bf Proof of \eqref{Abs_f}}

Now we establish \eqref{Abs_f}. There are two cases. The first case is $l_0(0)=0$. Then, applying \eqref{Abstract_a},
we obtain that
$$
\P\left( \sig=\infty \right)\geq\Pi_{k=0}^\infty\P\left( l_0(k+1)=0|l_0(k)=0 \right)\geq \Pi_{k=0}^\infty p_k.
$$
The second case is $l_0(0)=\infty$. Then
$$
\P\left( \sig=\infty \right)\geq\P\left( l_0(1)=1 \right)\Pi_{k=1}^\infty\P\left( l_0(k+1)=1|l_0(k)=1 \right).
$$
Since $\Hcal_0\leq R_0$, then applying \eqref{Abstract_a} and \eqref{Abstract_b}
$$
\P\left( \sig=\infty \right)\geq \Pi_{k=-1}^\infty p_k.
$$
Since $p_k>0$ and $1-p_k$ exponentially decreases, then the product converges and  in the two cases
\begin{equation}\label{Abs_l}
\P\left( \sig=\infty \right)\geq p_\infty =\Pi_{k=-1}^\infty p_k>0.
\end{equation}
Notice that \eqref{Abstract_a} implies
$$
\P\left( \sig = n \right)\leq \P\left( l_0(n+1)\not= n\,|\,l_0(n)=0\right)+ \P\left( l_0(n+1)\not= n\,|\,l_0(n)=1\right)
\leq 2c_1 e^{-\alpha_1 (n-1)T},
$$
which gives the first inequality of \eqref{Abs_e} and allows to conclude

\section{Properties of the CGL equation}

We are concerned with the stochastic Complex Ginzburg--Landau (CGL)
equations with Dirichlet boundary conditions:
\begin{equation}\label{Equation_Intro_glc_Dirichlet}
\Espace
\left \{
\begin{array}{rcl}
\glcDelta{u} & = & b(u)\frac{dW}{dt}+f , \\
u(t,x) & = & 0, \quad \mbox{   for   } x \in \delta D, \\
u(0,x) & = & u_0(x),
\end{array}
\right .
\end{equation}
where  $\eps >0$, $\eta >0$, $\lambda\in\{-1,1\}$ and where $D$ is an open  bounded set of $\R^d$ with sufficiently
regular boundary or $D=[0,1]^d$.
 Also $f$ is the deterministic part of the forcing term.
For simplicity in the redaction, we consider the case $f=0$. The generalisation to a square integrable $f$
 is easy.
 We say that it is the defocusing or the focusing equation when $\lambda$ is equal to
$1$ or $-1$, respectively.

We set
$$
A = -\Delta, \quad  D(A) = H^1_0(D) \cap H^2(D).
$$
Now we can write problem (\ref{Equation_Intro_glc_Dirichlet}) in the form
\begin{eqnarray}
\glc{u} & = & b(u)\frac{dW}{dt} ,\label{Equation_glc_zeta} \\
u(0) & = & u_0,  \label{Equation_glc_zeta_conditions_initiales}
\end{eqnarray}
where $W$ is a cylindrical Wiener process of $L^2(D)$.

The aim of this section is to
prove some properties which will be used in Section 3 to build a coupling such that the assumptions of Theorem
\ref{Th_Theo}  are true.

\subsection{Notations and main result}
\

We consider $(e_n,\mu_n)_{n\in \N^*}$ the couples of eigenvalues and eigenvectors of $A$ ($Ae_n=\mu_n$) such that $(e_n)_n$
is an Hilbertian basis of $L^2(D)$ and such that $(\mu_n)_n$ is an increasing sequence. We denote by $P_N$ and $Q_N$
the orthogonal projection in  $L^2(D)$  on the space $Sp(e_k )_{1\leq n}$ and on its complementary,
 respectively.

The first condition is a condition on the smoothness of the noise and a condition ensuring existence and uniqueness of solutions.

We will sometimes consider the $L^2(D)$ sub-critical condition:

{\bf H1 } {\it We assume that $0 < \sig < \frac{2}{d}\wedge \frac{3}{2}$.
Moreover $u_0 \in L^2(D)$ and b is bounded Lipschitz
$$
b: L^2(D) \to \mathcal{L}_2 (L^2(D),H^2(D)).
$$
}

We also consider the $H^1(D)$ sub-critical condition  when the equation is defocusing.

{\bf H1' }{\it If $d \leq 2$ we assume that $\sig >0$. If $d>2$, we assume that $0 < \sig < \frac{2}{d-2}$.
Moreover $\lambda=1$, $u_0 \in H^1(D)$ and b is bounded Lipschitz
$$
b: L^2(D) \to \mathcal{L}_2 (L^2(D),H^2(D)).
$$
}

We set, for $s\leq 2$,
$$B_s=\sup_u\abs{b(u)}^2_{\mathcal{L}_2 (L^2(D),H^s(D))}.$$

The second assumption means that $b$ only depends on its low modes.

{\bf H2 }{\it There exists $N_1$ such
that
$$
b(u)=b(P_{N_1} u).
$$}

The third condition is a structure condition on $b$. It is a slight generalisation of the usual assumption that $b(u)$
is diagonal in the basis $(e_n)_n$.

{\bf H3 }{\it There exists $N\geq N_1$, such that for any $u$,
$$
P_N b(u) Q_N =0, \quad Q_N b(u) P_N =0.
$$
Moreover $ P_N b(u)P_N$ is  invertible on $P_N H$ and
$$
\sup_u \abs{(P_N b(u)P_N)^{-1}} < \infty.
$$
}

In this section, we define by $\abs{\cdot}$, $\abs{\cdot}_p$, $\norm{\cdot}$ and $\norm{\cdot}_s$ the norm of
$L^2(D)$, $L^p(D)$, $H^1(D)$ and $H^s(D)$.

The Lyapunov structures are defined by
$$
\Espace
\begin{array}{rcl}
\Hcal^{L^2}&=&\abs{\cdot}^2,\\
 \Hcal^{H^1}&=&\frac{1}{2}\norm{\cdot}^2+\frac{1}{2\sigma+2}\abssig{\cdot}.
\end{array}
$$

The energies are defined by
$$
E_u^{L^2}(t,T)=\abs{u(t)}^2+\eps \int_T^t\norm{u(s)}^2 ds,
$$
and
$$
\Espace
E_u^{H^1}(t,T)=
\left \{
\begin{array}{l}
\Hcal^{H^1}(u(t))+ \frac{\eps}{2} \int_T^t\norm{u(s)}_2^2 ds +  \frac{\eta}{2} \int_T^t \absquatre{u(s)}ds
\\+
(\eta+\eps)\int_T^t \int_D \abs{u(s,x)}^{2\sig}\abs{\nabla u(s,x)}^2 dx ds,
\end{array}
\right.
$$
When $T=0$, we simply write $E_u(t)=E_u(t,0)$.

The first case is the $L^2$--subcritical  focusing or defocusing CGL equation with initial condition  in $L^2(D)$:

{\bf Case 1:}
\begin{itemize}
\item  {\bf H1}, {\bf H2} and {\bf H3} hold,
\item $\lambda\in\{-1,1\},\quad H=L^2(D)$,
\item $\Hcal=\Hcal^{L^2}=\abs{\cdot}_{L^2(D)}^2,\quad E_u=E_u^{L^2}$.
\end{itemize}

The second case is the $H^1$--subcritical defocusing CGL equation with initial condition in $H^1(D)$.

{\bf Case 2:}
\begin{itemize}
\item  {\bf H1'}, {\bf H2} and {\bf H3} hold,
\item $\lambda= 1,\quad H=H^1(D)$,
\item $\Hcal=\Hcal^{H^1}=\frac{1}{2}\norm{\cdot}_{H^1(D)}^2+\frac{1}{2\sig+2}\abs{\cdot}_{L^{2\sig+2}(D)}^{2\sig+2},
\quad E_u=E_u^{H^1}$.
\end{itemize}

When it is not precised, the results stated are true in both cases. It is well known that we have existence and uniqueness
of the solutions in both cases and that the solutions are strongly Markov process. We denote by $(\Pcal_t)_{t\in \R^+}$ the
Markov transition semi-group associated to the solutions of \eqref{Equation_glc_zeta}.

The aim of this article is to establish the following result
\begin{Theorem}[MAIN THEOREM]\label{Th_MAIN}
There exists $N_0(B_2,\eta,\eps,\sig,D)$ such that if $N\geq N_0$, then in cases 1 and 2, there exists a unique stationnary
 probability measure $\nu$ of $(\Pcal_t)_{t\in \R^+}$ on $L^2(D)$. Moreover, $\nu$ satisfies
\begin{equation}\label{Eq_MAIN_a}
\int_{H} \norm{u}^2_{H^2(D)}d\nu(u) < \infty,
\end{equation}
and for any $s \in [0,2)$, there exists $C_s>0$ and $\alpha_s$  such that for any $\mu \in \Pcal(H)$
\begin{equation}\label{Eq_MAIN_b}
\abs{\Pcal^*_t\mu-\nu}_{Lip_b(H^s(D))}^*\leq C_s e^{-\alpha_s t}\left(1+ \int_{H} \abs{u}^2_{L^2(D)} d\mu(u) \right).
\end{equation}
Furthermore, if $(u,W)$ is a weak solution of \eqref{Equation_glc_zeta}, \eqref{Equation_glc_zeta_conditions_initiales},
with $u_0$ taking value in $L^2(D)$ then for any $f \in Lip_b(H^s(D))$
\begin{equation}\label{Eq_MAIN_h}
\abs{\E f(u(t)) -\int_H f(u) d\nu(u)}\leq C_s \abs{f}_{Lip_b(H^s(D))}  e^{-\alpha_s t}\left(1+ \E \abs{u_0}^2_{L^2(D)} \right).
\end{equation}
\end{Theorem}

\begin{Remark}\label{Remark_MAIN2}
In case 1, \eqref{Eq_MAIN_b} is equivalent to \eqref{Eq_MAIN_h}. But in case 2, the Markovian transition semi-group make
 sense only if $u_0$ is taking value in $H= H^1(D)$ because strong existence and weak uniqueness may cause problem when
$u_0\in L^2(D)$. Hence \eqref{Eq_MAIN_b} make sense only if $\mu \in \Pcal(H^1(D))$ which means that $u_0 \in H^1(D)$.
\end{Remark}

\begin{Remark}\label{Remark_MAIN}
Assume that $B_s<\infty$ for $s$ sufficiently high. Let $k$ be a positive integer such that
$$
k \leq 2\sig+2,\textrm{ if } \sig \not \in \N,\textrm{ and } k\in \N \textrm{ if } \sig \in \N.
$$
Applying Remark \ref{Remark_H2} below and adapting the proof of Theorem \ref{Th_MAIN}, we obtain that
\eqref{Eq_MAIN_a} can be replaced by
\begin{equation}\label{Eq_MAIN_c}
\int_{H} \norm{u}^2_{H^k(D)}d\nu(u) < \infty,
\end{equation}
and \eqref{Eq_MAIN_b} is true for any $s$  real number such that
$$
s < [2\sig+2],\textrm{ if } \sig \not \in \N,\textrm{ and } s\in \R \textrm{ if } \sig \in \N,
$$
where $[\cdot]$ denote the integer part.

The condition on $k$ and $s$ comes from the lack of derivability of the non-linear part of the CGL equation. Assume that we
replace $\abs{u}^{2\sig}u$ by $g(\abs{u}^2)u$ where
\begin{itemize}
\item $g$ is infinitely continuously differientiable,
\item $g(x)=x^\sig $ for $x\geq x_0$,
\item $g$ is increasing and $g(0)=0$.
\end{itemize}
Hence Theorem \ref{Th_MAIN}, \eqref{Eq_MAIN_c} and \eqref{Eq_MAIN_b} are true for any $k$ and $s$.
\end{Remark}
\subsection{Properties of the solutions}
\

In this subsection, we state some properties proved in the next subsections. These are used in Section 3 to
 apply Theorem \ref{Th_Theo} in order to establish Theorem \ref{Th_MAIN}.

First, we recall the following result.
\begin{Proposition}\label{Th_iso}
In the two previous cases, there exists a mesurable map
$$
\Phi :  C((0,T);P_N H)\times C((0,T);Q_N H^{\frac{d+1}{2}}(D)) \times  H     \to  C((0,T);Q_N H),
$$
such that for any $(u,W)$ solution of \eqref{Equation_glc_zeta} and \eqref{Equation_glc_zeta_conditions_initiales}
$$
Q_N u   =  \Phi( P_N u ,Q_N W,u_0)\quad \textrm{ on } [0,T].
$$
Moreover $\Phi$ is a non-anticipative functions of $(P_N u ,Q_N W)$.
\end{Proposition}
Proposition \ref{Th_iso} can be proved by applying a fix point argument and by taking into account that the limit of
a sequence of measurable maps is measurable.

We have the so-called Foias-Prodi estimates.
\begin{Proposition}[Foias-Prodi estimate]\label{Prop_Foias_Prodi_Hun}
Let $u_1$ and $u_2$ be two solutions of the CGL system
\eqref{Equation_glc_zeta} associated with  Wiener process $W_1$ and $W_2$
 respectively.
 If
 \begin{equation}\label{Hyp_Prop_Foias_Prodi_Hun}
 P_N u_1(t) = P_N u_2(t), \quad Q_N W_1(t)=Q_N W_2(t),
 \mbox{  for  }  T_0 \leq t \leq T,
 \end{equation}
where $N$ is a non-negative integer, then
 \begin{equation}\label{Resultat_Prop_Foias_Prodi_Hun}
\abs{r(t)}_H \leq \abs{r(T_0)}_H
\exp \left (  -\frac{\eps \mu_{N+1}}{2} (t-T_0) +
              c_1 \sum_{i=1}^2 E_{u_i}(t,T_0)   \right ),
 \end{equation}
where $r=u_1-u_2$ and $  T_0 \leq t \leq T$ and where $c_1>0$ only depends on $\eps$, $\eta$, $\sig$, $D$.
\end{Proposition}
We deduce immediately a very usefull Corollary.
\begin{Corollary}\label{Cor_Foias_Hun}
For any $B$, there exists $N_0'(B,\eta,\eps,D,\sig)$
 such that under the assumptions of
 Proposition \ref{Prop_Foias_Prodi_Hun}, under the assumption $N \geq N_0'$ and under the assumption
$$
E_{u_i}(t,T_0)\leq  \rho +B (t-T_0),\quad i=1,2
$$
we obtain that
$$
\abs{r(t)}_H \leq \abs{r(T_0)}_H
\exp \left (  - 2(t-T_0)
+c_1\rho  \right ).
$$
where $c_1$ is the constant of Proposition \ref{Prop_Foias_Prodi_Hun}.
\end{Corollary}
Then, by proving analogous result to the previous Corollary, we obtain the Drift estimate which, in Section 3, will ensures the Novikov
condition and will allow to apply the Girsanov Formula.
\begin{Lemma}[Drift estimate]\label{Prop_Novikov_Hun}
 For any $B$, there exists $N_0"(B,\eta,\eps,D,\sig)$ such that for any
$u_1,u_2$ solutions of the CGL system
\eqref{Equation_glc_zeta} associated with $W_1$ and $W_2$ and  for any   $N>N_0"$
 \begin{equation}\label{Resultat_Prop_Novikov_Hun}
\int_{T_0}^\tau\abs{P_N(\B{u_1(s)}-\B{u_2(s)})}^2ds
\leq K_N\abs{r(0)}^2e^{c\rho-3T_0},
 \end{equation}
 where $T>T_0\geq  0$ and $\rho, C,\alpha>0$,
 where  $K_N,c$ only depend on $B$, $C$, $\alpha$, $\eps$, $\eta$, $\sig$, $D$, $N$
 and where we have denoted by $\tau$ the value
$$
\tau=T_0\vee \inf \left( t\in [0,T] \left|
\begin{array}{l}
 E_{u_1}(t)\geq  \rho +B t \mbox{ or } E_{u_2}(t)\geq  \rho +C\left( 1+t^\alpha\right) \mbox{ or }\\
 P_N u_1(t) \not = P_N u_2(t)\mbox{ or  }
Q_N W_1(t) \not =Q_N W_2(t)
\end{array}
 \right. \right).$$
\end{Lemma}
Now we set
$$
N_0=N_0'\vee N_0".
$$
In order to apply the previous Lemmas and Corollary, we establish the two following results.
\begin{Proposition}[Exponential estimate for the growth of solution]\label{Prop_majoration_energie_sous_lineaire_Hun}
Assume that $u$ is a solution of
\eqref{Equation_glc_zeta}, \eqref{Equation_glc_zeta_conditions_initiales}
associated with a Wiener process $W$.
Then, for any $0 \leq T_0 < T \leq \infty$
$$
\P \left (  \sup_{t \in [T_0,T[} \left (   E_u(t) - B t  \right )
\geq  \Hcal(u_0)+\rho  \right) \leq e^{- \gamma_0 \rho -  3T_0},
$$
where $B$ only depends on $B_2$, $\sigma$, $\eta$, $\eps$.
\end{Proposition}
\begin{Proposition}\label{Prop_majoration_energie_sous_lineaire_Hunbis}
Assume that $u$ is a solution of
\eqref{Equation_glc_zeta}, \eqref{Equation_glc_zeta_conditions_initiales}
associated with a Wiener process $W$. For any $u_0^2$, we define $\widetilde u$ by
$$
\widetilde u=P_N u+\phi\left(P_N u,Q_N W,u_0^2\right).
$$
Then, there exists $\alpha\geq 1$ such that for any $N$, there exists $C_N$,
$$
\P \left (  \sup_{t \in [0,T[} \left (   E_{\widetilde u}(t) - C_Nt^\alpha  \right )
\geq  C_N\left(1+\Hcal(u_0)+\Hcal(u_0^2)^\alpha+\rho \right) \right) \leq 2e^{- \gamma_0 \rho},
$$
 for any $0 \leq  T \leq \infty$ and any $u_0^2$.
\end{Proposition}

Let $u_1$ and $u_2$ be two solutions of \eqref{Equation_glc_zeta} that
correspond to deterministic intial value $u_0^1$ and $u_0^2$, respectively.
\begin{Lemma}[The Lyapunov structure]\label{lem_Lyapounov_norm}
 There exists $\alpha>0$ and $C_k>0$ such that for any $k$
$$
 \E \mathcal H(u_i(t))^k \leq \mathcal H(u_0^i)^k e^{-\alpha k t } + \frac{C_k}{2},
 $$
and for any stopping time $\tau$
$$
\E \mathcal H(u_i(\tau))^k 1_{\tau < \infty}\leq  \Hcal(u_0^i)^k +
C_k\left(1+\E \left(\tau 1_{\tau < \infty}\right)\right).
$$
\end{Lemma}
Using Lemma \ref{lem_Lyapounov_norm} and  Chebyshev's inequality, we obtain
\begin{Lemma}\label{lem_Lyapounov_proba}
 If $R_0 \geq  (\Hcal(u_0^1)+\Hcal(u_0^2) )\vee C_1 $, then
$$
\P\left ( \mathcal H(u_1(t))+\Hcal(u_2(t)) \geq  4 C_1  \right ) \leq \frac{1}{2},
$$
providing
$t \geq  \theta_1(R_0)= \frac{1}{\alpha} \ln \frac{R_0}{C_1}$.
\end{Lemma}

Then, in  the second case, we control $\Hcal(u(t))$ by $\abs{u_0}^2$.
\begin{Proposition}\label{Prop_Ctrl}
It is assumed that $u$ is a solution of
\eqref{Equation_glc_zeta}, \eqref{Equation_glc_zeta_conditions_initiales}
associated with a Wiener process $W$.
Then, for any $T >0$
$$
\E \Hcal(u(T)) \leq A + B T +\frac{C}{T} \abs{u_0}^2 ,
$$
where $A$, $B$ and $C$  only depends on $B_2$, $\sigma$, $\eta$, $\eps$.
\end{Proposition}
Now, we claim that in the two cases, we can control the norm of solutions in Sobolev spaces by the norm in $L^2$.

\begin{Proposition}\label{Prop_H2}
Let $k$ be a positive integer less than $2$.
There exist $\gamma_k>1$ only depending on $k$, $\sig$ and $d$ and $C_k>0$ and $c_k>0$ only depending on
$k$, $(B_s)_s$, $\sig$, $d$, $\eps$ and $\eta$ such that for any $T>0$ and $t>0$
$$
\E \left( \norm{u(T+t)}_k^2 + \int_T^{T+t}\norm{u(s)}_{k+1}^2 ds \right)^{\frac{2}{\gamma_k}}
\leq c_k\frac{1}{T}\abs{u_0}^2+C_k(1+T+t).
$$
\end{Proposition}
Hence, applying a Chebyshev inequality, we obtain
\begin{Corollary}\label{Cor_H2}
Let $k$ be a positive integer less than $2$ and $\delta>0$. There exist $\gamma>0$  only depending on
$k$, $\sig$ and $d$ and $C_\delta>0$  only depending on
$\delta$, $k$, $(B_s)_s$, $\sig$, $d$, $\eps$ and $\eta$ such that for any  $t>0$
$$
\P\left( \norm{u(t)}_k \geq e^{\delta t} \right)\leq C_\delta e^{-\frac{\delta}{\gamma}t}\left( \abs{u_0}^2+1 \right)
$$
\end{Corollary}

\begin{Remark}\label{Remark_H2}
Assume that $B_s<\infty$ for $s$ sufficiently high. The proof of Proposition \ref{Prop_H2} can be adapted to $k$ a positive integer such that
$$
k \leq 2\sig+2,\textrm{ if } \sig \not \in \N,\textrm{ and } k\in \N \textrm{ if } \sig \in \N,
$$
and then Corollary \ref{Cor_H2} is true for such a $k$.

 The condition on $k$ comes from the fact that
$\abs{\cdot}^{\sig}$ is not $C^\infty$ on $0$. As in Remark \ref{Remark_MAIN}, if we replace
$\abs{\cdot}^{\sig}$ by a nice function which coincides with $\abs{\cdot}^{\sig}$ on $[x_0,\infty)$, we can
establish those results for any $k$.
\end{Remark}

\subsection{Foias-Prodi and Drift estimates}
\

The proofs in the first case are closely related to the proofs in the second case, but are simpler.
 That is the reason why
 we only give the proof in the second case.

{\bf Proof of Proposition \ref{Prop_Foias_Prodi_Hun} in the second case.}

We denote $u_1-u_2$ by $r$.

{\it Step 1.} This step is devoted to the proof of
\begin{equation}\label{Preuve_Prop_Foias_Prodi_a_Hun}
I=(  (\eta + \i)(\B{u_2}-\B{u_1})  , A r ) \leq
\frac{\eps}{2} \norm{r}_2^2  + c \norm{r}^2
\sum_i\absquatre{u_i}.
\end{equation}
We recall the following estimate
\begin{equation}\label{AC}
\abs{\B{x}-\B{y}}
\leq c\abs{x-y}(\abs{x}^{2\sig}+\abs{y}^{2\sig}).
\end{equation}

Applying  H\"older inequality  and then  \eqref{AC} gives
$$
I\leq \norm{r}_2 \abs{\B{u_2}-\B{u_1}}_2 \leq \norm{r}_2 \abs{\left(\sum_{i=1}^2\abs{u_i}^{2\sig}\right) r}_2.
$$
Let $s \in (1,2)$ such that $\frac{4\sig}{2-s}=4\sig+2$. Applying  once more H\"older inequality and then
the Sobolev embedding $H^{s}(D) \subset L^{4\sig+2}(D)$ gives
$$
I \leq \norm{r}_2\abs{r}_{4\sig+2}\sum_{i=1}^2\abs{u_i}_{4\sig+2}^{2\sig}
\leq \norm{r}_2\norm{r}_{s}\sum_{i=1}^2\abs{u_i}_{4\sig+2}^{2\sig},
$$
which yields by the interpolatory inequality $\norm{.}_s\leq\norm{.}_2^{s-1}\norm{.}^{2-s}$ and then an
arithmetic-geometric inequality
$$
I \leq \norm{r}_2^s\norm{r}^{2-s}\sum_i\abs{u_i}_{4\sig+2}^{2\sig}\leq \frac{\eps}{2} \norm{r}_2^2  + c \norm{r}^2
\sum_i\absquatre{u_i}.
$$

{\it Step 2.} We now establish \eqref{Resultat_Prop_Foias_Prodi_Hun}.

 Taking into
account \eqref{Hyp_Prop_Foias_Prodi_Hun}, we see that $r$ satisfies
the equation
\begin{equation}\label{Preuve_Prop_Foias_Prodi_b_Hun}
\chal{r} = (\eta + \i) Q_N (\B{u_2}-\B{u_1}).
\end{equation}
Taking the scalar product of \eqref{Preuve_Prop_Foias_Prodi_b_Hun} by
$-2A r$, we obtain:
 \begin{equation}\label{Preuve_Prop_Foias_Prodi_c_Hun}
\chalcarreH{r}{2} = 2(  (\eta + \i)(\B{u_2}-\B{u_1})  , Ar ).
\end{equation}
Taking into account \eqref{Preuve_Prop_Foias_Prodi_a_Hun},
\eqref{Preuve_Prop_Foias_Prodi_c_Hun} gives :
 \begin{equation}\label{Preuve_Prop_Foias_Prodi_d_Hun}
\chalcarreH{r}{} \leq c \norm{r}^2
\sum_i\absquatre{u_i} .
\end{equation}
Since $r \in Q_N H$, then  $\mu_{N+1}\norm{r}^2 \leq \norm{r}_2^2$ and it follows from
\eqref{Preuve_Prop_Foias_Prodi_d_Hun} that
 \begin{equation}\label{Preuve_Prop_Foias_Prodi_e_Hun}
\frac{d\norm{r}}{dt} + \eps \mu_{N+1} \norm{r}^2 \leq c \norm{r}^2
\sum_i\absquatre{u_i} .
\end{equation}
Applying Gromwall Lemma to \eqref{Preuve_Prop_Foias_Prodi_e_Hun}, we obtain
\eqref{Resultat_Prop_Foias_Prodi_Hun}.
\carre

{\bf Proof of Lemma \ref{Prop_Novikov_Hun} in the second case.}

\noindent We first state the following Lemma which strengthen Proposition \ref{Prop_Foias_Prodi_Hun}.
\begin{Lemma}\label{Prop_Foias_Prodi_Hunbis}
Let $u_1$ and $u_2$ be two solutions of the CGL system
\eqref{Equation_glc_zeta} associated with  $W_1$ and $W_2$
 respectively.
 If
 \begin{equation}\label{Hyp_Prop_Foias_Prodi_Hunbis}
 P_N u_1(s) = P_N u_2(s), \quad Q_N W_1(s)=Q_N W_2(s),
 \mbox{  for  any }  s \in (T_0,t),
 \end{equation}
where $N$ is a non-negative integer, then
 \begin{equation}\label{Resultat_Prop_Foias_Prodi_Hunbis}
\abs{r(t)}_{L^2} \leq \abs{r(0)}_{L^2}
\exp \left (  -\frac{\eps \mu_{N+1}}{2} t +
              c_1  E_{u_1}(t)   \right ),
 \end{equation}
where $r=u_1-u_2$  and where $c_1>0$ only depends on $\eps$, $\eta$, $\sig$, $D$.
Moreover,  for any $B$, there exists $N_0"(B,\eta,\eps,D,\sig)$ such that $N\geq N_0"$ and
 \begin{equation}\label{Hyp_Prop_Foias_Prodi_Hunter}
E_{u_1}(t)\leq  \rho +B t
 \end{equation}
imply
 \begin{equation}\label{Hyp_Prop_Foias_Prodi_Hun4}
\abs{r(t)}_{L^2} \leq \abs{r(0)}_{L^2}
\exp \left (  - 2t
+c_1\rho  \right ),
 \end{equation}
where $c_1$ is the constant of Proposition \ref{Prop_Foias_Prodi_Hun}.
\end{Lemma}
For the first case,  this result is Proposition 1.1.6 of \cite{ODASSO0}. For the second case the proof is the same.

\noindent {\bf Sketch of the proof of Lemma \ref{Prop_Foias_Prodi_Hunbis}.}

\noindent The proof  of Lemma \ref{Prop_Foias_Prodi_Hunbis} is similar to the proof of Proposition
\ref{Prop_Foias_Prodi_Hun}. Indeed it is sufficient to prove
\begin{equation}\label{Rqbis}
I'=-(  (\eta + \lambda i)(\B{u_2}-\B{u_1})  , r ) \leq c \abs{
\abs{u_1}^{ 2\sigma}  \abs{ r}^2 }_1 .
 \end{equation}
 to establish Lemma \ref{Prop_Foias_Prodi_Hunbis}. We prove \eqref{Rqbis} as follows. Remarking that
$$
\B{u_2}-\B{u_1}=\abs{u_2}^{2\sigma} r+u_1(\abs{
u_2}^{2\sigma}-\abs{u_1}^{2\sigma}) ,
$$
and
$$
\abs{u_1(\abs{
u_2}^{2\sigma}-\abs{u_1}^{2\sigma})}\leq c'\abs{u_1}(\abs{
u_2}^{2\sigma-1}+\abs{u_1}^{2\sigma-1})\abs{r},
$$
we obtain
$$
I'\leq-\eta\abs{\abs{u_2}^{2\sigma}  \abs{ r}^2 }_1+
c'\abs{\abs{u_1}^{ 2\sigma}  \abs{ r}^2 }_1
+c'\abs{u_1\abs{u_2}^{2\sigma-1} r^2}_1.
$$
Applying arithmetico-geometric inequality to the last term of the previous equality,
 we obtain for $\sig\geq\frac{1}{2}$
$$
c'\abs{u_1\abs{u_2}^{2\sigma-1} r^2}_1\leq
\eta \abs{\abs{u_2}^{2\sigma}  \abs{ r}^2 }_1 +c" \abs{\abs{u_1}^{2\sigma}  \abs{ r}^2 }_1.
$$
We infer \eqref{Rqbis} for $\sig\geq\frac{1}{2}$ from the two previous inequalities.

\noindent To obtain  \eqref{Rqbis} when
$\sig<\frac{1}{2}$, one remark that $D$ is the union of $\{x|\abs{u_1(x)}\geq\abs{u_2(x)}\}$ and
$\{x|\abs{u_1(x)}<\abs{u_2(x)}\}$.
 Treating the first set is trivial. The treatement done before works for the second set.
\carre

Let us set
$$I=\int_{T_0}^\tau\abs{P_N(\B{u_1(s)}-\B{u_2(s)})}^2ds .$$
Applying Lemma \ref{Prop_Foias_Prodi_Hunbis}
 with the same $N_0"$, we obtain
\begin{equation}\label{Preuve_Novikov_a_Hun}
\abs{r(t)} \leq \abs{r(0)}
\exp \left (  - 2t
+c_1\rho  \right ),  \mbox{  for  }    \tau\geq t \geq 0.
 \end{equation}
Noticing that, since we work in a finite dimensional space, all the norm are equivalent. Hence
 there exists $K_N$ such that
\begin{equation}\label{Preuve_Novikov_b_Hun}
I \leq K_N \int_{T_0}^\tau\abs{\B{u_1(s)}-\B{u_2(s)}}_1^2ds .
 \end{equation}
It follows from \eqref{AC} and H\"older inequality that
\begin{eqnarray}
\abs{\B{u_1(s)}-\B{u_2(s)}}_1^2   &\leq & c\abs{\left(\sum_{i=1}^2\abs{u_i(s)}^{2\sig}\right)\abs{r(s)}}_1^2,\nonumber\\
 & \leq & c\left(\sum_{i=1}^2\abs{u_i(s)}_{4\sig}^{4\sig}\right)\abs{r(s)}^2, \nonumber
\end{eqnarray}
which yields, by applying an
arithmetico-geometric inequality,
\begin{equation}\label{Preuve_Novikov_c_Hun}
\abs{\B{u_1(s)}-\B{u_2(s)}}_1^2  \leq  c\left(1+\sum_{i=1}^2\abs{u_i(s)}_{4\sig+2}^{4\sig+2}\right)\abs{r(s)}^2.
 \end{equation}
Combining \eqref{Preuve_Novikov_a_Hun}, \eqref{Preuve_Novikov_b_Hun}  and \eqref{Preuve_Novikov_c_Hun}
 and then an integration by parts, we obtain
\begin{eqnarray}
I &\leq &K_N \abs{r(0)}^2 \int_{T_0}^\tau  \exp \left (  - 4t
+c_1\rho  \right ) \left(1+\sum_{i=1}^2\abs{u_i(s)}_{4\sig+2}^{4\sig+2}\right) ds \nonumber, \\
& \leq & K_N \abs{r(0)}^2 \int_{T_0}^\tau  \exp \left (  - 4t
+c_1\rho  \right )\left(1+ \sum_{i=1}^2\int_{T_0}^t\abs{u_i(s)}_{4\sig+2}^{4\sig+2}ds\right) dt \nonumber,\\
& \leq & K_N \abs{r(0)}^2 \int_{T_0}^\tau  \exp \left (  - 4t
+c_1\rho  \right )\left(1+ 2\rho + B t+C(1+t^\alpha)\right) dt \nonumber,\\
& \leq & K_N \abs{r(0)}^2 \int_{T_0}^\tau  \exp \left (  - 3t
+2 c_1\rho  \right ) dt \nonumber,
\end{eqnarray}
which allows us to conclude.
\carre
\subsection{An exponential estimate for the growth of solution}
\

As in the previous subsection, we only give the proof of  Propositions \ref{Prop_majoration_energie_sous_lineaire_Hun}
  in the second case.

We set
$$
E'_u(t)=
\left \{
\begin{array}{l}
\frac{1}{2}\norm{u(t)}^2+\frac{1}{2\sigma+2}\abssig{u(t)}+ \eps \int_0^t\norm{u(s)}_2^2 ds +  \eta \int_0^t \absquatre{u(s)}ds
\\+\\
(\eta+\eps)\int_0^t \int_D (1+\chi(u\nabla \bar u))\abs{u(s,x)}^{2\sig}\abs{\nabla u(s,x)}^2 dx ds,
\end{array}
\right.
$$
where $\chi(z)=2\sigma \Re e \left( \frac{\Re e z}{z} \right)$.
Applying Ito's Formula to $\Hcal(u)=\frac{1}{2}\norm{u}^2+\frac{1}{2\sigma+2}\abssig{u}$, we obtain
\begin{equation} \label{Preuve_nrj_sslin_a_Hun}
 E'_u(t)=\Hcal(u_0)+ M_1(t)+M_2(t)+I_1(t)+I_2(t),
\end{equation}
where we have denoted
$$
\begin{array}{ll}
M_1(t)= \int_0^t (-\Delta u(s) , b(u(s)) dW(s)), & M_2(t)= \int_0^t (\B{ u(s)} , b(u(s))) dW(s)),\\
I_1(t)=\frac{1}{2} \int_0^t \abs{b(u(s))}^2_{\mathcal{L}^2(L^2(D),H^1(D))} ds, &   I_2(t)=\frac{1}{2} \int_0^t \sum_{i=1}^2\abs{g_i(u(s))}^2_{\mathcal{L}^2(L^2(D))} ds ,
\end{array}
$$
where
$$
g_i(u)(k)=f_i(u)(b(u)h)\quad f_1(u)(k)=\abs{u}^{\sig}\times k,\quad f_2(u)(k)=\sqrt{2\sig}\abs{u}^{\sig-1}\Re e (\bar u \times k ).
$$
H\"older estimate and Sobolev Embedding give
$$
\sum_{i=1}^2\abs{f_i(u)}^2_{\mathcal{L}(H^1(D,L^2(D)))} \leq c\abs{u}_{4\sig+2}^{4\sig},
$$
which yields
$$
\sum_{i=1}^2\abs{g_i(u)}^2_{\mathcal{L}^2(L^2(D))} \leq c\abs{u}_{4\sig+2}^{4\sig}B_1,
$$
and thus by an arithmetico-geometric inequality
\begin{equation}\label{Eq_exp_Hun_maj_Ideux}
I_2(t)\leq c B_1 t + \frac{\eta}{4} \int_0^t \absquatre{u}ds
\end{equation}
Notice that
$$
<M_1>(t)=\int_0^t \abs{b(u(s)^*A u(s)}^2 ds,
$$
which gives
\begin{equation}\label{Eq_exp_Hun_maj_Mun}
<M_1>(t)\leq B_0 \int_0^t \norm{u(s)}_2^2 ds.
\end{equation}
Moreover
$$
<M_2>(t)=\int_0^t \abs{b(u(s))^* \B{u(s)}}^2 ds.
$$
Since
$$\abs{b(u(s))^* \B{u(s)}}^2\leq B_0 \abs{ \B{u(s)}}^2\leq c B_0\absquatre{u},$$
we obtain
\begin{equation}\label{Eq_exp_Hun_maj_Mdeux}
<M_2>(t)\leq c B_0 \int_0^t \absquatre{u}ds .
\end{equation}
Noticing that $<M_1+M_2>\leq 2(<M_1>+<M_2>)$, $I_1(t)\leq B_1 t$ and $\chi(z)\geq 0$ for any $z \in \C$, it follows from \eqref{Preuve_nrj_sslin_a_Hun}, \eqref{Eq_exp_Hun_maj_Ideux}, \eqref{Eq_exp_Hun_maj_Mun} and \eqref{Eq_exp_Hun_maj_Mdeux} that
\begin{equation}\label{Eq_exp_Hun_maj_nrj}
E_u(t)-\Hcal(u_0)-B t \leq M(t)-\frac{\gamma_0}{2}<M>(t),
\end{equation}
where $M=M_1+M_2$, $B'=c(B_0+B_1)$ and $\gamma_0=\frac{\eta\vee \eps}{8B_0(1+c)}$.
Thus
$$
\P(\sup_{t\in \R^+}\left (   E_u(t) - B t  \right )\geq \Hcal(u_0)+\rho' ) \leq e^{-\gamma_0 \rho'}\E e^{\gamma_0 M(t)-\frac{\gamma_0^2}{2}<M>(t)}\leq e^{-\gamma_0 \rho'},
$$
which allows to conclude by  setting $\rho'=\rho+3\frac{T_0}{\gamma_0}$ and $B'=B+\frac{3}{\gamma_0}$.

We do not give the proof of Proposition \ref{Prop_majoration_energie_sous_lineaire_Hunbis} because it is easilly
deduced from the proof of Proposition \ref{Prop_majoration_energie_sous_lineaire_Hun}. Actually, Ito Formulas
 associated to a solution $u$ are also true if we replace $u$ by $\widetilde u$ and $b(P_Nu)dW$ by
$b(P_N \widetilde u)dW+ P_N(\B u-\B{\widetilde u}) dt$.
Hence to establish Proposition \ref{Prop_majoration_energie_sous_lineaire_Hunbis}, it is sufficient to bound
the additionnal term by using the equivalence of the norms in finite-dimensionnal spaces and by applying
 Proposition \ref{Prop_majoration_energie_sous_lineaire_Hunbis} to bound terms containing $u$.

\subsection{ The Lyapunov structure}
\

Now, we prove  Lemma \ref{lem_Lyapounov_norm} in the second case.
Using the computation of the energy previously done, we obtain that there exixts $C_1$ such that
$$
d\mathcal H(u_i(t))+ \frac{\eps}{2}  \norm{u_i(t)}_2^2  dt
+ \frac{\eta}{4\sigma+2}  \abs{u_i(t)}_{4\sigma+2}^{4\sigma+2} dt
\leq dM+ C_1 d t
$$
Applying Ito Formula to $\Hcal(u_i)^k$ and  controlling $d<M>$ as above by $\norm{u_i(t)}_2^2  dt$ and
$\abs{u_i(t)}_{4\sigma+2}^{4\sigma+2} dt $, we obtain that there exists $\alpha_0$ such that
\begin{equation}\label{Eq_Lyap_c}
\Espace
\begin{array}{ll}
d\mathcal H(u_i(t))^k+ \alpha_0 k \Hcal(u_i)^{k-1}&\left(\norm{u_i(t)}_2^2+
 \abs{u_i(t)}_{4\sigma+2}^{4\sigma+2}\right) dt \\
&\quad \quad \quad \leq k \Hcal(u_i(t))^{k-1} dM+ C_k d t.
\end{array}
\end{equation}
Taking into account that $\mu_1\norm{.}^2 \leq \norm{.}_2^2$ and that there exist $\beta>0$ such that $\beta \abssig{.}
\leq \norm{.}_2^2+\abs{.}_{4\sigma+2}^{4\sigma+2}$, we obtain that there exists $\alpha>0$ such that
\begin{equation}\label{Eq_Lyap_a}
d\mathcal H(u_i(t))^k+ \alpha k \Hcal(u_i)^k dt
\leq k \Hcal(u_i(t))^{k-1} dM+ C_k d t,
\end{equation}
which  yields, by integrating and taking the expectation, the second inequality of Lemma \ref{lem_Lyapounov_norm}.

Now, applying \eqref{Eq_Lyap_a}, we obtain that
\begin{equation}\label{Eq_Lyap_b}
\mathcal H(u_i(t))^k\leq  \Hcal(u_0^i)^k e^{-\alpha k t}+k \int_0^t e^{-\alpha k(t-s)} \Hcal(u_i(s))^{k-1} dM(s)+ C_k .
\end{equation}
which  yields, by taking the expectation, the first inequality of Lemma \ref{lem_Lyapounov_norm}.

\subsection{Control of $\Pcal_T \Hcal$ by $\abs{.}^2$ in  the second case}
\

Now, we prove Proposition \ref{Prop_Ctrl}.
Taking the expectation on \eqref{Eq_exp_Hun_maj_nrj}, we obtain that for any $T>t>0$
$$
\E \Hcal(u(T)) \leq \E \Hcal(u(t)) +B (T-t).
$$
Integrating over $[0,T]$ gives
\begin{equation}\label{Eq_Ctrl}
\E \Hcal(u(T)) \leq \frac{1}{T} \E\int_0^T \Hcal(u(t))dt +B T.
\end{equation}
Applying Ito Formula to $\abs{u}^2$ and taking the expectation, we obtain
$$
\E \abs{u(t)}^2+2\eps \int_0^t\E\norm{u(s)}^2 ds + 2 \eta \int_0^t \E\abssig{u(s)}ds =
 \abs{u_0}^2 +  \int_0^t \E\abs{b(u(s))}^2_{\mathcal{L}_2 (L^2(D)) } ds.
$$
 Applying {\bf H1'}, we obtain
$$
\E\int_0^T \Hcal(u(t))dt \leq C\abs{u_0}^2+ AT,
$$
and by \ref{Eq_Ctrl}
$$
\E \Hcal(u(T)) \leq A + B T +\frac{C}{T} \abs{u_0}^2 .
$$

\subsection{$H^1$ and $H^2$ estimates}
\

We first establish that
\begin{equation}\label{Eq2.7.0}
\E\norm{u(T)}^{2}+\eps\int_0^T\E \norm{u(s)}_2^{2}ds\leq \norm{u_0}^{2}+c_1\abs{u_0}^{\alpha_1}+B_1' T,
\end{equation}
and that
\begin{equation}\label{H1_e}
\E\norm{u(T)}^2\leq c\left(1+\frac{1}{T}\abs{u_0}^2+\abs{u_0}^{2k}+T\right).
\end{equation}
In the second part of the proof, we establish that there exists $\gamma_0 >0$ such that
\begin{equation}\label{Eq2.7.2}
\E\norm{u(t)}^2_2+\eps\int_0^t\E \norm{u(s)}_3^2ds\leq \norm{u_0}_2 ^2+ c\norm{u_0}^{\gamma_0}+C( t+1).
\end{equation}
We deduce from H\"older inequality that
\begin{equation}\label{H2_f}
\E\left(\norm{u(t)}^2_2+\eps\int_0^t \norm{u(s)}_3^2ds\right)^\frac{2}{\gamma_0}\leq c\norm{u_0}_2^{2}+C( t+1).
\end{equation}
and
\begin{equation}\label{Eq2.7.3}
\E\norm{u(T)}_2^2\leq c\left(1+\frac{1}{T}\norm{u_0}^2+T\right).
\end{equation}
Hence, combining \eqref{H1_e}, \eqref{H2_f} and \eqref{Eq2.7.3}, we obtain
\begin{equation}\label{H2_g}
\E\left(\norm{u(T+t)}^2_2+\eps\int_T^{T+t} \norm{u(s)}_3^2ds\right)^\frac{2}{\gamma_0}\leq
c\frac{1}{T}\abs{u_0}^{2}+\abs{u_0}^{2k}+C( T+ t+1).
\end{equation}
Applying H\"older inequality allows to conclude.

{\bf Proof of \eqref{Eq2.7.0} and  \eqref{H1_e}}

Note that \eqref{Eq2.7.0} and \eqref{H1_e} have already been demonstrated in the second case. Then it remains to establish \eqref{Eq2.7.0}
 in the first case, when $\lambda=-1$.

Remark that Ito's Formula applied to $\abs{u}^{2k}$ gives
\begin{equation}\label{L2_a}
\E\left( \abs{u(t)}^{2k}+\eta k \int_0^t \abs{u(s)}^{2(k-1)}\abs{u}_{2\sig+2}^{2\sig+2}ds \right)\leq \abs{u_0}^{2k} +B_k" t.
\end{equation}

Taking the scalar product between \eqref{Equation_glc_zeta} and $2(-\Delta)u$ gives
\begin{equation}\label{H1_a}
d\norm{u}^2+2\eps\norm{u}_2^2dt\leq 2((-\Delta u),b(u)dW)+2(\Delta u,(\eta+\lambda \i)\abs{u}^{2\sig}u)dt+ B_1 dt.
\end{equation}
We deduce from Schwartz inequality that
$$
2(\Delta u,(\eta+\lambda \i)\abs{u}^{2\sig}u)\leq c\norm{u}_2 \abs{u}_{4\sig+2}^{2\sig+1}.
$$
The Gagliardo-Niremberg inequality gives
$$
2(\Delta u,(\eta+\lambda \i)\abs{u}^{2\sig}u)\leq c\norm{u}_2^{1+\frac{\sig d}{2}} \abs{u}^{2\sig+1-\frac{\sig d}{2}}.
$$
Finally, since $\sig d<2$, then we can deduce from a arithmetico-geometric inequality that
$$
2(\Delta u,(\eta+\lambda \i)\abs{u}^{2\sig}u)\leq \eps\norm{u}_2^2+c \abs{u}^{2\frac{4\sig+2-\sig d}{2-\sig d}}.
$$
We infer from \eqref{H1_a} that
$$
d\norm{u}^2+\eps\norm{u}_2^2dt\leq 2((-\Delta u),b(u)dW)+c \abs{u}^{2\frac{4\sig+2-\sig d}{2-\sig d}}dt+B_1dt,
$$
and then
$$
\E\norm{u(t)}^2+\eps\int_0^t\E \norm{u(s)}_2^2ds\leq \norm{u_0}^2
+   c \int_0^t \E \abs{u(s)}^{2\frac{4\sig+2-\sig d}{2-\sig d}}ds+B_1t.
$$
Applying \eqref{L2_a}, we obtain for a well-chosen $k'$
\begin{equation}\label{H1_d}
\E\norm{u(t)}^2+\eps\int_0^t\E \norm{u(s)}_2^2ds\leq c\left(\norm{u_0}^2+\abs{u_0}^{2k'}+ T\right).
\end{equation}
Using the same argument as in the last subsection gives \eqref{H1_e}.

{\bf Proof of \eqref{Eq2.7.2}, \eqref{H2_f} and \eqref{Eq2.7.3}}

Taking the scalar product between \eqref{Equation_glc_zeta} and $2(-\Delta)^2 u$ gives
\begin{equation}\label{H2_a}
d\norm{u}_2^2+2\eps\norm{u}_3^2dt\leq 2((-\Delta u)^2,b(u)dW)-2((-\Delta)^2 u,(\eta+\lambda \i)\abs{u}^{2\sig}u)dt+B_2dt.
\end{equation}
We deduce from an integration by part and Schwartz inequality that
\begin{equation}\label{H2_b}
-2((-\Delta)^2 u,(\eta+\lambda \i)\abs{u}^{2\sig}u)\leq c\norm{u}_3 \abs{\nabla \left(u \abs{u}^{2\sig}\right)}.
\end{equation}
H\"older inequality gives
$$
\abs{\nabla \left(  u \abs{u}^{2\sig}\right)}\leq \abs{\nabla u}_p \abs{u}_{2\sig q}^{2\sig},
$$
where $\frac{1}{p}+\frac{1}{q}=\frac{1}{2}$.
We choose $s$, $p$ and $q$ such that
$$
\frac{1}{p}=\frac{1}{2}-\frac{s}{d}, \quad \frac{1}{2\sig q}=0\vee \left( \frac{1}{2}-\frac{1}{d}\right).
$$
Since $\sig \leq \frac{2}{d-2}$, then $s\in [0,2)$.
 Hence the Sobolev embeddings $H^s(D)\to L^p(D)$ and $H^1(D)\to L^{2\sig q}(D)$ imply
$$
\abs{\nabla u \abs{u}^{2\sig}}\leq \norm{ u}_{1+s} \norm{u}^{2\sig},
$$
Then, we deduce from \eqref{H2_b}, an interpolatory inequality that
$$
-2((-\Delta)^2 u,(\eta+\lambda \i)\abs{u}^{2\sig}u)
\leq c\norm{u}_3^{1+\frac{s}{2}} \norm{u}^{2\sig+1-\frac{s}{2}}.
$$
An arithmetico-geometric inequality gives
\begin{equation}\label{H2_d}
-2((-\Delta)^2 u,(\eta+\lambda \i)\abs{u}^{2\sig}u)
\leq \eps\norm{u}_3^2 +c \norm{u}^{\beta},
\end{equation}
with  $\beta>0$.
We infer from \eqref{H2_a} and \eqref{H2_d} that
\begin{equation}\label{H2_e}
d\norm{u}_2^2+\eps\norm{u}_3^2dt\leq 2((-\Delta u)^2,b(u)dW)+c \norm{u}^{\beta}dt+B_2dt.
\end{equation}
Hence, we deduce \eqref{Eq2.7.2} from \eqref{H2_e}.
Then, applying H\"older inequality, we obtain \eqref{H2_f}.
Using the same argument as in the last subsection gives \eqref{Eq2.7.3}.

\section{The coupling of CGL}

Recall that, as in the last section, we consider the two cases devellopped in subsection 2.1 and use the properties
stated in subsection 2.2. In this section, we make an other assumption

{\bf H4} \quad \quad \quad \quad \quad $ N\geq N_0$,

where $N_0$ has been defined after Corollary \ref{Cor_Foias_Hun} and  Lemma \ref{Prop_Novikov_Hun}.

In this section, we apply Theorem \ref{Th_Theo}. Then we obtain there exists a unique invariant probability measure on
$H$ and that there exists $c>0$ and $\alpha>0$
\begin{equation}\label{Cgl_a}
\P\left( \abs{u_1(t)-u_2(t)}_{H}>c e^{-\alpha t}  \right)\leq c e^{-\alpha t}\left( 1+\Hcal(u_0^1)+\Hcal(u_0^2) \right).
\end{equation}
 Recalling Corollary \ref{Cor_H2}, we obtain for any $\delta>0$,
\begin{equation}\label{Cgl_b}
\P\left( \norm{u_i(t)}_{H^2(D)} \geq e^{\delta t} \right)\leq C_{\delta} e^{-\frac{\delta}{\gamma}t}
\left( \abs{u_0^i}_{L^2(D)}^2+1 \right)
\end{equation}
Combining \eqref{Cgl_a}, \eqref{Cgl_b} and using an interpolatory inequality between $L^2(D)$ and $H^2(D)$, we obtain
that for any $s\in [0,2)$, there exists $\alpha_s>0$ and $C_s>0$ such that
$$
\P\left( \norm{u_1(t)-u_2(t)}_{H^s(D)}>c e^{-\alpha_s t}  \right)\leq C_s e^{-\alpha_s t}
 \left( 1+\sum_{i=1}^2\left(\abs{u_0^i}_{L^2(D)}^{2}+\Hcal(u_0^i)\right)\right),
$$
which implies
$$
\norm{\Pcal_t^* \mu -\nu}_{Lip_b(H^s(D))}^*\leq C_s e^{-\alpha_s t}
 \left( 1+\int_{H}\left(\abs{u}_{L^2(D)}^{2}+\Hcal(u)\right)d\mu(u)\right).
$$
Now it remains to conclude the second case,
we consider $(u,W)$ a weak solution and we apply Proposition \ref{Prop_Ctrl}
$$
\E\left( \Hcal(u(T))+\abs{u(T)}_{L^2(D)}^{2}\right)\leq \frac{1}{T} \E\abs{u_0}^{2}_{L^2(D)}+C(1+T).
$$
which implies for all cases
$$
\abs{\E f(u(t)) -\int_H f(u) d\nu(u)}
\leq c_s \abs{f}_{Lip_b(H^s(D))} e^{-\alpha_s t}\left( 1+\E\abs{u_0}^{2}_{L^2(D)} \right),$$
 for any  $s< 2$, for any $f \in Lip_b(H^s(D))$.

It follows from this discussion that it suffices to prove that Theorem \ref{Th_Theo} can be applied and that \eqref{Cgl_a}
 holds. Then Theorem \ref{Th_MAIN} is proved.

\subsection{Preliminaries}
\

We set
$
\abs{\cdot}=\abs{\cdot}_H
$
 and
$$
X=P_N u,\quad Y=Q_N u,\quad \beta=P_N W,\quad \eta= Q_N W,\quad \sig_l=P_N b P_N, \quad \sig_h=Q_N b Q_N,
$$
and
$$
\begin{array}{rcl}
f(X,Y)&=&(\eta+\lambda \i)P_N\left( \abs{X+Y}^{2\sig}(X+Y) \right),\\
g(X,Y)&=&(\eta+\lambda \i)Q_N\left( \abs{X+Y}^{2\sig}(X+Y) \right).
\end{array}
$$
Now, taking into account {\bf H2} and {\bf H3}, the system has the form
\begin{equation}\label{CGL_Eq_un_neuf}
\left \{
\Espace \begin{array}{lcl}
d X+ (\eps+\i)A X dt +f(X,Y)dt & = & \sig_l(X)d\beta , \\
d Y+ (\eps+\i)A Y dt +g(X,Y)dt & = & \sig_h(X)d\eta  , \\
\lefteqn{ X(0)=x_0 , \quad Y(0)=y_0.}
\end{array}
\right .
\end{equation}
Recall that {\bf H3} states that
\begin{equation}\label{CGL_Eq_un_neuf_ter}
 \textrm{There exists } \sig_0>0 \textrm{ such that,} \; \abs{\left(\sig_l(x)\right)^{-1}}\leq \frac{1}{\sig_0}, \;
 \textrm{ for any } x\in P_N H.
\end{equation}
Now we can define $l_0$
$$
l_0(k)=\min\left\{l \in \{0,...,k\} | P_{l,k}\right\},
$$
 where $\min \phi =\infty$ and
$$
(P_{l,k})
\left \{
\Espace
\begin{array}{l}
X_1(t)=X_2(t), \quad \eta_1(t)=\eta_2(t), \quad \forall\; t \in [lT,kT],\\
\Hcal_l\leq d_0 ,\quad i=1,2,\\
E_{u_i}(t+lT,lT)\leq \aleph 1_{t<T}+B t+1_{i=2}1_{t\leq T}C_N(1+t^\alpha), \quad \forall\; t \in [0,(k-l)T],
\end{array}
\right.
$$
where $B,\alpha,C_N$ are defined in Propositions \ref{Prop_majoration_energie_sous_lineaire_Hun} and
\ref{Prop_majoration_energie_sous_lineaire_Hunbis},
where $\aleph$ will be chosen later and where
$$
\Hcal_k=\Hcal(u_1(kT))+\Hcal(u_2(kT)).
$$
 Notice that \eqref{Abstract_ater} is obvious.
 Corollary \ref{Cor_Foias_Hun} and {\bf H4} gives
\begin{equation}\label{CGL_Eq_abis}
l_0(k)=l \quad \textrm{implies} \quad   \abs{u_1(t)-u_2(t)}\leq C(d_0) e^{-(t-lT)}, \; \textrm{ for any } t \in [lT,kT],
\end{equation}
 and we have establish \eqref{Abstract_abis}.
Lemma \ref{lem_Lyapounov_norm} implies the Lyapunov structure \eqref{Abstract_c}.

 From now on we say that $(X_1,X_2)$ are coupled at $kT$ if $l_0(k)\leq k$, in other words
 if $l_0(k)\not = \infty$.
Now it remains to build a coupling such that \eqref{CGL_Eq_a} and \eqref{CGL_Eq_b} holds, where
\begin{equation}\label{CGL_Eq_a}
\left\{
\Espace
\begin{array}{l}
\forall \; d_0,\;\exists  \;  p_0(d_0)>0, \; (p_i)_{i \in \N^*}, \; T_0(d_0)>0  \; \textrm{ such that for any } l\leq k , \\
\P\left(l_0(k+1)=l \; | \; l_0(k)=l \right) \geq p_{k-l}, \textrm{ for any }  T\geq T_0(d_0), \\
1-p_i \leq  e^{-i T}, \; i \in \N^*,
\end{array}
\right.
\end{equation}
 and, for any $(R_{0},d_0)$ sufficiently large,
\begin{equation}\label{CGL_Eq_b}
\left\{
\Espace
\begin{array}{l}
\exists \; T^*(R_0) >0  \textrm{ and } p_{-1}>0 \textrm{ such that for any }T\geq T^*(R_0)\\
\P\left(l_0(k+1)=k+1 \; | \; l_0(k)=\infty ,\; \mathcal H_k \leq R_0 \right)
 \geq p_{-1},
\end{array}
\right.
\end{equation}

These properties imply \eqref{Abstract_a} and \eqref{Abstract_b} and  Theorem \ref{Th_Theo} can be applied.

As in the example of section 1.2, we remark that by induction, it suffices to construct a probability space $(\Omega_0,\F_0,\P_0)$ and two measurable couples of functions
 $(\omega_0,u_0^1,u_0^2)\to(V_i(\cdot,u_0^1,u_0^2))_{i=1,2}$ and  $(V'_i(\cdot,u_0^1,u_0^2))_{i=1,2}$ and  such that, for any
$(u_0^1,u_0^2)$, $(V_i(\cdot,u_0^1,u_0^2))_{i=1,2}$ and $(V'_i(\cdot,u_0^1,u_0^2))_{i=1,2}$ are two couplings of
$(\Dr(u(\cdot,u_0^i),W))_{i=1,2}$ on $[0,T]$. Indeed, we first set
$$
u_i(0)=u_0^i,\quad W_i(0)=0, \quad i=1,2.
$$
Assuming that we have built $(u_i,W_i)_{i=1,2}$ on $[0,kT]$, then we take $(V_i)_i$ and $(V'_i)_i$ as above independant
 of $(u_i,W_i)_{i=1,2}$ on $[0,kT]$ and set
\begin{equation}\label{CGL_Eq_m_bis}
\Espace
\left(u_i(kT+t),W_i(kT+t)\right) =\left \{
\Espace \begin{array}{ll}
V_i(t,u_1(kT),u_2(kT)) & \textrm{ if } l_0(k) \leq k,\\
V'_i(t,u_1(kT),u_2(kT))& \textrm{ if } l_0(k) =\infty,
\end{array}
\right.
\end{equation}
for any $t \in [0,T]$.

\subsection{Proof of \eqref{CGL_Eq_a}}
\

The essential difference between this proof and the proof of \eqref{Eq_a} in the example in section 1.2 is that a cut-off
 is used to control the energy.

To build $(V_i(\cdot,u_0^1,u_0^2))_{i=1,2}$, we apply Proposition \ref{Prop_Matt} to
$$
\Espace
\begin{array}{l}
E=C((0,T); H)\times C((0,T);H^{-\frac{d}{2}-1}(D)),\\
F=C((0,T);P_N H)\times C((0,T);Q_N H^{-\frac{d}{2}-1}(D)),\\
f_0\left(u,W\right)= (X,\eta), \\
\mu_i=\Dr(u(\cdot,u_0^i),W), \quad \textrm{ on } [0,T].
\end{array}
$$
Remark that if we set $\nu_i=f_0^*\mu_i$, we obtain
$$
\nu_i=\Dr(X(\cdot,u_0^i),\eta), \quad \textrm{ on } [0,T].
$$
We set
$$
(Z_i,\xi_i)=f_0(V_i), \quad i=1,2.
$$

Then $(V_i(\cdot,u_0^1,u_0^2))_{i=1,2}$ is a coupling of $(\mu_1,\mu_2)$ such that
$((Z_i,\xi_i)(\cdot,u_0^1,u_0^2))_{i=1,2}$ is a maximal coupling of $(\nu_1,\nu_2)$.

 We first use a Girsanov formula to estimate $I_p$, where
$$
\Espace
\begin{array}{lcl}
I_p&=&\int_{A_{k,l}} \left ( \frac{d\nu_2}{d \nu_1} \right )^{p+1} d\nu_2,\\
A_{k,l} &= & \{(Z,\xi)\,|\,\tau_{k,l}=T\},\\
\tau_{k,l}&=& \inf\left\{t\in[0,T]\;|\;  E_{\hat u_i}(t+kT,lT)>\aleph 1_{k=l}+B(t+(k-l)T)\right.\\
&&\quad\quad\quad\quad\quad\quad\quad\quad\quad\quad\quad\quad\quad\quad
\left.+1_{i=2}1_{k=l}C_N(1+t^\alpha),\;i\in\{1,2\}\right\},
\end{array}
$$
where
$$
\BLANC{\Espace
\begin{array}{l}
 \hat E^{k,l}_{Z,\xi,u_0,u}(t)=E_{\hat u}(kT+t,lT),\\}
 \hat u_i = u_i \; \textrm{  on } [0,kT], \quad \hat u_i(kT+\cdot)=Z+\Phi(Z,\xi,u_0^i)\; \textrm{  on } [0,T].
\BLANC{\end{array}}
$$
Then, using Lemma \ref{lem_norm_var_sens}, we establish \eqref{CGL_Eq_a}.

We consider a couple of  $(u_i,W_i)_{i=1,2}$,
two solutions of \eqref{CGL_Eq_un_neuf} on $[0,kT]$ and a trajectory of $(u_i,W_i)_{i=1,2}$
 such that $l_0(k)=l$.  We set
$$
x=X_1(kT)=X_2(kT), \quad y_i=Y_i(kT), \quad i=1,2.
$$
Let $W=(\beta,\xi)$ a cylindrical Wiener process
 defined on a probability space $(\Omega,\F,\P)$. We denote by $Z$ the unique solution of the truncated equation
\begin{equation}\label{CGL_Eq_i_un}
\Espace
\left\{
\begin{array}{lcl}
d Z+ (\eps+\i)A Z dt +1_{t\leq\tau_{k,l}}f(Z,\Phi(Z,\xi,(x,y_1)))dt  =  \sig_l(Z) d\beta , \\
\lefteqn{ Z(0)=x.}
\end{array}
\right.
\end{equation}
  We denote by $\lambda_1$  the distribution of
$(Z,\xi)$ under the probability $\P$.

\noindent We set $\widetilde{\beta}(t)=\beta(t)+\int_0^t d(s) dt $ where
$$
 d(t)=1_{t\leq\tau_{k,l}}\left(\sig_l(Z(t))\right)^{-1}\left(f(Z(t),\Phi(Z,\xi,(x,y_2))(t))-f(Z(t),\Phi(Z,\xi,(x,y_1))(t))\right).
$$
Then $Z$ is a solution of
\begin{equation}\label{CGL_Eq_i_deux}
\Espace
\left\{
\begin{array}{l}
d Z+ (\eps+\i)A Z dt +1_{t\leq\tau_{k,l}}f(Z,\Phi(Z,\xi,(x,y_2)))dt  =  \sig_l(Z) d\widetilde{\beta} , \\
\lefteqn{ Z(0)=x.}
\end{array}
\right.
\end{equation}
The drift estimate in Lemma \ref{Prop_Novikov_Hun} ensures  that
\begin{equation}\label{CGL_Eq_i_trois}
\int_0^T\abs{d(t)}^2 dt\leq c d_0 \sig_0^{-2}  \exp\left( -3(k-l)T+c\aleph 1_{k=l} \right).
\end{equation}
Hence the Novikov condition is satisfied and  the Girsanov formula
can be applied. Then we set
$$
d\widetilde{\P}= \exp\left(\int_0^T d(s) dW(s)-\frac{1}{2} \int_0^T \abs{d(s)}^2 dt  \right)d\P
$$
 We deduce from the Girsanov formula that $\widetilde{\P}$ is a probability under which $(\widetilde{\beta},\xi)$ is a
 cylindrical Wiener process and we denote by $\lambda_2$ the law of $(Z,\xi)$ under $\widetilde{\P}$.
 Moreover, remarking that
\begin{equation}\label{CGL_Eq_jbis}
\lambda_i (A_{k,l}\cap\cdot) =\nu_i(A_{k,l}\cap\cdot), \quad i=1,2,
\end{equation}
we obtain
\begin{equation}\label{CGL_Eq_jter}
I_p\leq I_p'\leq \E
\exp\left( c_p \int_0^T \abs{d(s)}^2 dt \right),
\end{equation}
where
$$
I_p'=\int_F \left ( \frac{d\lambda_2}{d \lambda_1} \right )^{p+1} d\lambda_2,\\
$$
Then it follows from \eqref{CGL_Eq_i_trois} that
\begin{equation}\label{CGL_Eq_k}
I_p\leq I_p' \leq
\exp\left( c_p \sig_0^{-2}d_0  e^{ -3(k-l)T +c\aleph 1_{k=l}} \right).
\end{equation}
Notice that
$$
\norm{\lambda_1-\lambda_2}_{var}=\int_F \abs{ \frac{d \lambda_2}{d \lambda_1}  -1}d\lambda_2
\leq \sqrt{\int \left ( \frac{d \lambda_2}{d \lambda_1} \right )^{2} d\lambda_2-1}.
$$
We infer from \eqref{CGL_Eq_k} that, for $T\geq T_3(d_0)=2\ln \left(c_p \sig_0^{-2}d_0\right)$,
$$
\norm{\lambda_1-\lambda_2}_{var}\leq \frac{1}{2}e^{ -2(k-l)T}.
$$
Using \eqref{CGL_Eq_jbis}, we obtain for $k>l$
$$
\norm{\nu_1-\nu_2}_{var}\leq \norm{\lambda_1-\lambda_2}_{var}+\sum_{i=1}^2\nu_i({A_{k,l}^i})\leq
 \frac{1}{2}e^{ -2(k-l)T}+\sum_{i=1}^2\nu_i({A_{k,l}^i}).
$$
where
$$
A_{k,l}^i=\left\{(Z,\xi)\,\left|\, E_{Z+\phi(Z,\xi,u_0^i)}(t,lT)\leq B(t+(k-l)T) \; \textrm{ for any } t\in [0,T]\right.  \right\}.
$$
 Applying Lemma \ref{lem_norm_var_sens} to the maximal coupling $(Z_1,Z_2)_{i=1,2}$ of $(\nu_1,\nu_2)$ gives for $k>l$
\begin{equation}\label{CGL_Eq_l}
\P\left( (Z_1,\xi_1) \not = (Z_2,\xi_2) \right) \leq \norm{\nu_1-\nu_2}_{var}\leq \frac{1}{2}e^{ -2(k-l)T}+
\sum_{i=1}^2\nu_i({A_{k,l}^i}).
\end{equation}
Using \eqref{CGL_Eq_m_bis} and \eqref{CGL_Eq_l}, we obtain that  on $l_0(k)=l$
$$
\P\left( (X_1,\eta_1) \not = (X_2,\eta_2) \textrm{ on } [kT,(k+1)T]  \,\left| \, \F_{kT} \right.\right)
\leq  \frac{1}{2}e^{ -2(k-l)T}+2\P({B_{l,k}}|\F_{kT}),
$$
where
$$
B_{l,k}=\left\{ E_{u_i}(t,lT)\leq B(t-lT), \; \textrm{ for any } t\in [kT,(k+1)T],\; i=1,2  \right\}.
$$
Noticing that for $k>l$
$$
\{l_0(k+1)=l\}=\{l_0(k)=l\}\cap\{(X_1,\eta_1)  = (X_2,\eta_2) \textrm{ on } [kT,(k+1)T]\}\cap B_{l,k}.
$$
and integrating over $l_0(k)=l$ gives for $T\geq T_1(d_0)$ and for $k>l$
$$
\P\left( l_0(k+1)\not=l \,|\,l_0(k)=l\right)
\leq \frac{1}{2}e^{ -2(k-l)T}+3\P({B_{l,k}}\,|\,l_0(k)=l),
$$
and then
$$
\P\left( l_0(k+1)\not=l,\; l_0(k)=l \,|\,l_0(l)=l\right)
\leq \frac{1}{2}e^{ -2(k-l)T}+3\P({B_{l,k}}\,|\,l_0(l)=l).
$$
The exponential estimate for growth of the solution (Proposition \ref{Prop_majoration_energie_sous_lineaire_Hun})
 gives that for $T$ sufficiently high
\begin{equation}\label{CGL_Eq_e}
\P\left( l_0(k+1)\not=l,\; l_0(k)=l \,|\,l_0(l)=l\right)
\leq \exp(-2(k-l)T).
\end{equation}

Now, it remains to consider the case $k=l$, we apply Lemmas  \ref{lem_norm_var_sens} and
 \ref{lem_tech_coup_inf} to
 $(Z_i,\xi_i)_{i=1,2}$ which gives
$$
\P\left( (Z_1,\xi_1)=(Z_2,\xi_2), \; A_{l,l}^2 \right)\geq \left( \nu_1\wedge\nu_2 \right)(A_{l,l})\geq
 \left(1-\frac{1}{p}\right) (pI_p)^{-\frac{1}{p-1}}\nu_1(A_{l,l})^{\frac{p}{p-1}}.
$$
Choosing $\aleph$ sufficiently high and applying the exponential for growth of the solution
 (Propositions \ref{Prop_majoration_energie_sous_lineaire_Hun} and
\ref{Prop_majoration_energie_sous_lineaire_Hunbis}), we obtain
$$
\nu_1(A_{l,l})\geq\frac{1}{2},
$$
and then applying \eqref{CGL_Eq_jter} and fixing $p>1$,
$$
\P\left( (Z_1,\xi_1)=(Z_2,\xi_2), \; A_{l,l} \right)\geq p_0(d_0)>0.
$$
That gives
\begin{equation}\label{CGL_Eq_ebis}
\P\left( l_0(l+1)=l\,|\,l_0(l)=l\right)\geq p_0(d_0)>0.
\end{equation}
Since
$$
\P\left( l_0(k)\not=l|l_0(l)=l \right) \leq \sum_{n=l}^{k-1} \P\left( l_0(n+1)\not=l,\; l_0(n)=l    \,|\,l_0(l)=l\right),
$$
then, by applying \eqref{CGL_Eq_e} and \eqref{CGL_Eq_ebis}, we obtain
$$
\P\left( l_0(k)\not=l|l_0(l)=l \right) \leq 1-p_0+ \sum_{n=1}^\infty \exp(- 2nT)\leq 1-p_0 +
\frac{\exp(- 2T)}{1-\exp(- 2T)},
$$
which implies that for $T\geq T_0(d_0)$
\begin{equation}\label{CGL_Eq_eter}
\P\left( l_0(k)=l|l_0(l)=l \right) \geq \frac{p_0}{2},
\end{equation}
 Combining \eqref{CGL_Eq_e}, \eqref{CGL_Eq_ebis} and \eqref{CGL_Eq_eter}, we establish \eqref{CGL_Eq_a} for $T$
sufficiently high.

\subsection{Proof of \eqref{CGL_Eq_b}}
\

As in the example of Section 1.2, The Lyapunov structure gives that it is sufficient to find $d_0>0$, $\widetilde p>0$,
  $R_1>4 K_1$ and
  a coupling $( V_i(\cdot,u_0^1,u_0^2))_{i=1,2}$ of
$( \mu_1, \mu_2)$, where
$$
 \mu_i=\Dr (u(\cdot,u_0^i),W), \; \textrm{ on } [0,1], \quad i= 1,2,
$$
and such that
\begin{equation}\label{CGL_Eq_b_b}
\P\left( Z_1(1,u_0^1,u_0^2)= Z_2(1,u_0^1,u_0^2),\;\sum_{i=1}^2\Hcal( u_i(1,u_0^1,u_0^2))\leq d_0
 \right)\geq \widetilde p,
\end{equation}
where
$$
 V_i(\cdot,u_0^1,u_0^2)=\left(  u_i(\cdot,u_0^1,u_0^2),  W_i(\cdot,u_0^1,u_0^2) \right),
\;  u_i(\cdot,u_0^1,u_0^2)=\left(\begin{array}{c}  Z_i\\ G_i \end{array}\right), \; i=1,2.
$$

Now we fix $R_1>4K_1$ and consider a cimetery value $\Delta$ (some people prefer calling it a heaven value).
To build $( V_i(\cdot,u_0^1,u_0^2))_{i=1,2}$, we apply Proposition \ref{Prop_Matt} to
$$
\Espace
\begin{array}{l}
E=C((0,1); H)\times C((0,1);H^{-\frac{d}{2}-1}(D)),\\
F=\left(P_N H\times C((0,1);Q_NH^{-\frac{d}{2}-1}(D))\right)\cup\{\Delta\},\\
f_0\left(u,W\right)= X(1) 1_{A}(X,\eta)+ \Delta 1_{A^c}(X,\eta),
\end{array}
$$
and to $ \mu_i$ where
$$
\Espace
\begin{array}{lcl}
A &= & \{(X,\eta)\,|\,\tau=1\},\\
\tau&=& \inf\left\{t\in[0,1]\;|\;  E_{X+\Phi(X,\eta,u_0^i)}(t)>\aleph+B t+1_{i=2}C_N(1+t^\alpha),\;i\in\{1,2\}\right\}.
\end{array}
$$
We set $\nu_i=f_0^*\mu_i$.
Then $( V_i(\cdot,u_0^1,u_0^2))_{i=1,2}$ is a coupling of $(\mu_1,\mu_2)$ such that
$(Z_i(1,u_0^1,u_0^2))_{i=1,2}$ is a maximal coupling of $(\nu_1,\nu_2)$.

 Now, we define
$$
f_1\left(u,W\right)= (X,\eta)\quad \textrm{ and }\quad
f_2\left(X,\eta\right)= X(1) 1_{A}(X,\eta)+ \Delta 1_{A^c}(X,\eta),
$$
and we set $\theta_i=f_1^* \mu_i$ for $i=1,2$.
Now we consider $(\hat \theta_1,\hat \theta_2)$ such that
$ \theta_i(A\cap\cdot)$ is equivalent to $\hat \theta_i(A\cap\cdot)$ for $i=1,2$ and such that
$(\hat \nu_1,\hat \nu_2)=(f_2^*\hat \theta_1,f_2^*\hat \theta_2)$ are two equivalent measures.  Then by applying
 two Schwartz inequalities, we obtain that
\begin{equation}\label{CGL_Eq_b_a}
I_p \leq \left( J_{2p+2}^1\right)^{\frac{1}{2}}
\left(J_{4p}^2\right)^{\frac{1}{4}}\left(\hat I_{4p+2}\right)^{\frac{1}{4}},
\end{equation}
where
$$
\Espace
\begin{array}{rclrcl}
I_p &=& \int_{B'} \left ( \frac{d \nu_1}{d \nu_2} \right )^{p+1} d\nu_2,&
J_p^1 &=& \int_A \left ( \frac{d \theta_1}{d \hat\theta_1} \right )^{p} d\hat\theta_1,\\
\hat I_p &=& \int_{B'} \left ( \frac{d \hat\nu_1}{d \hat\nu_2} \right )^{p} d\hat\nu_2,&
J_p^2 &=& \int_A \left ( \frac{d \hat\theta_2}{d \theta_2} \right )^{p} d\hat\theta_2,
\end{array}
$$

\noindent Let us consider $\bar Z_i$ the unique solution of
\begin{equation}\label{CGL_Eq_b_d_un}
\Espace
\left\{
\begin{array}{lcl}
d \bar Z_i+ (\eps+\i)A \bar Z_i dt +1_{t\leq\tau}f(\bar Z_i,\Phi(\bar Z_i(\cdot),\xi(\cdot),u_0^i))dt & = & \sig_l(\bar Z_i) d\beta_i , \\
\lefteqn{ \bar Z_i(0)=x_0^i.}
\end{array}
\right.
\end{equation}
Taking into account \eqref{CGL_Eq_i_un}, we denote by $\lambda_i$ the distribution of
$(\bar Z_i,\xi_i)$ under the probability $\P$ and  we obtain
\begin{equation}\label{CGL_Eq_b_d_i_huit}
 \theta_i(A\cap\cdot)=\lambda_i(A\cap\cdot).
\end{equation}
\noindent We set $\widetilde{\beta_i}(t)=\beta_i(t)+\int_0^t d_i(s) dt $ where
\begin{equation}\label{CGL_Eq_b_d_j}
 d_i(t)=-1_{t\leq\tau}(\sig_l(\bar Z_i(t)))^{-1}f(\bar Z_i(t),\Phi(\bar Z_i(\cdot),\xi(\cdot),u_0^i)(t)).
\end{equation}
Then $\bar Z_i$ is a solution of
\begin{equation}\label{CGL_Eq_b_d_i_deux}
\Espace
\left\{
\begin{array}{lcl}
d \bar Z_i+ (\eps+\i)A \bar Z_i dt  & = & \sig_l(\bar Z_i) d\widetilde{\beta_i} , \\
\lefteqn{ \bar Z_i(0)=x_0^i.}
\end{array}
\right.
\end{equation}
Since the energy is bounded and $\sig_l$ is bounded below, then $d$ is uniformly bounded. Hence, the Novikov condition is satisfied and the Girsanov formula
can be applied. Then we set
$$
d\widetilde{\P_i}= \exp\left(\int_0^T d_i(s) dW(s)-\frac{1}{2} \int_0^T \abs{d(s)}^2 dt  \right)d\P
$$
 We deduce from the Girsanov formula that $\widetilde{\P}$ is a probability under which $(\widetilde{\beta},\xi)$ is a
 cylindrical Wiener process. We denote by $\hat \theta_i$ the law of $(\bar Z_i,\xi_i)$ under $\widetilde{\P_i}$.
 Moreover using \eqref{CGL_Eq_b_d_i_huit}, we obtain
\begin{equation}\label{CGL_Eq_b_d_jbis}
J_p^1\vee J_p^2\leq \E
\exp\left( c_p \int_0^T \abs{d(s)}^2 dt \right)\leq C(p,\aleph,R_1).
\end{equation}
We set $\hat \nu_i=f_2^*\hat \theta_i$ for $i=1,2$. It is classical that $\hat \nu_i$
has a density $q(x_0^i,z)$ with respect to lebesgue measure $dz$, that $q$ is continuous for the couple
$(x_0^i,z)$, where $x_0^i$ is the initial value and where $z$ is the target value and that $q>0$. Then, we can
 bound $q$ and $q^{-1}$ uniformly on
$\Hcal(x_0^i)\leq R_1$ and on $z\in B' =\left\{\Hcal(z)\leq C\right\}$ provided $C=C(\aleph)$.
 It allows us to bound $\hat I_p$ and then $I_p$. Actually, 
$d_1\geq d_1(\aleph)$  implies
\begin{equation}\label{CGL_Eq_b_e0}
A\subset B, 
\end{equation}
where
$$
B=\left\{(Z,\xi)\left|\Hcal(Z(1)+\phi(Z,\xi,u_0^i)(1))\leq d_1,\,i=1,2\right.\right\}.
$$
Hence it follows that for $d_1\geq d_1(\aleph)$
\begin{equation}\label{CGL_Eq_b_e}
I_p\leq C'(p,\aleph,R_1)<\infty.
\end{equation}
Now we apply Lemma  \ref{lem_tech_coup_inf} and  \ref{lem_norm_var_sens}:
\begin{equation}\label{CGL_Eq_b_ebis}
\P\left( Z_1(1)= Z_2(1),\;  (A\cap B')^2\right)\geq
\left(1-\frac{1}{p}\right)(pI_p)^{-\frac{1}{p-1}}\nu_1(B')^{\frac{p}{p-1}}.
\end{equation}

\noindent
 We deduce from Propositions \ref{Prop_majoration_energie_sous_lineaire_Hun} and
 \ref{Prop_majoration_energie_sous_lineaire_Hunbis} and from $C(\aleph)\to\infty$ when $\aleph\to\infty$
  that $\aleph$ sufficiently high gives
\begin{equation}\label{CGL_Eq_b_eter}
\nu_1(B')\geq \frac{1}{2}.
\end{equation}
Combining \eqref{CGL_Eq_b_e0}, \eqref{CGL_Eq_b_e}, \eqref{CGL_Eq_b_ebis} and \eqref{CGL_Eq_b_eter} gives
 for $d_1\geq d_1(\aleph)$
\begin{equation}\label{CGL_Eq_b_f}
\P\left( Z_1(1)=Z_2(1),\;  B^2\right)\geq \widetilde p= \widetilde p(p,\aleph,R_1)>0.
\end{equation}
\BLANC{Notice that
\begin{equation}\label{CGL_Eq_b_g}
\Espace
\begin{array}{l}
\P\left( Z_1(1)= Z_2(1),\;  \Hcal( u_i(1))\leq d_1,\; i=1,2\right)
\geq\\
\quad\quad\quad\quad\quad\quad\quad\quad\quad\quad\quad\quad\quad\P\left( (Z_1(1),\eta_1)= (Z_2(1),\eta_2),\;  B\right).
\end{array}
\end{equation}
Using the Lyapunov structure, we obtain that
\begin{equation}\label{CGL_Eq_b_h}
\P\left(\Hcal( u_i(T_1))\geq d_2\right)\leq \frac{R_1+K_1}{d_1+d_2}.
\end{equation}}
Taking into account the definition of $\phi$ and choosing 
 $d_0=2d_1$, 
 it follows that \eqref{CGL_Eq_b_b} holds.

\appendix

\section{ Proof of Lemma \ref{lem_norm_var_sens}}
Let $(Y_i)_i$ be a coupling of $(\mu_i)_i$. Let $\Gamma$ be a measurable set. There exists $(\Gamma_i)_i$ such that
$$
\Gamma = \bigcup_i \Gamma_i, \quad \bigcap_i \Gamma_i = \emptyset \quad (\mu_2-\mu_1)^+(\Gamma_2)=0, \quad (\mu_1-\mu_2)^+(\Gamma_1)=0.
$$
It follows from $a\wedge b=a-(a-b)^+$ and $(\mu_1-\mu_2)^+(\Gamma_1)=0$ that
$$
(\mu_1 \wedge \mu_2)(\Gamma_1)=\mu_1(\Gamma_1)-(\mu_1-\mu_2)^+(\Gamma_1)=\mu_1(\Gamma_1)= \P(Y_1 \in \Gamma_1).
$$
Symetricly, we obtain
$(\mu_1 \wedge \mu_2)(\Gamma_2)= \P(Y_2 \in \Gamma_2)$.

Thus, it follows from $\Gamma = \bigcup_i \Gamma_i$ and $ \bigcap_i \Gamma_i = \emptyset$ that
$$
(\mu_1 \wedge \mu_2)(\Gamma)= \P(Y_1 \in \Gamma_1)+\P( Y_2 \in \Gamma_2)\geq \sum_{i=1}^2\P(Y_1=Y_2, \;Y_1 \in \Gamma_i).
$$
Since  $\Gamma = \bigcup_i \Gamma_i$ and $ \bigcap_i \Gamma_i = \emptyset$ , then
\begin{equation}\label{Prf_norm_var_sens_a}
(\mu_1 \wedge \mu_2)(\Gamma)\geq \P (Y_1=Y_2, Y_1 \in \Gamma).
\end{equation}
Then it follows from $\norm{\mu_1-\mu_2}_{var}=1-(\mu_1 \wedge \mu_2)(E)$ that
$$
\norm{\mu_1-\mu_2}_{var} \leq \P(Y_1\not = Y_2).
$$
We have equality only if \eqref{Prf_norm_var_sens_a} appears for $\Gamma=E$, which is true only if \eqref{Prf_norm_var_sens_a} appears for any $\Gamma$.
For any measure $\mu$ on $(E,\mathcal{E})$, we denote by $\mathbf{\mu}$  the measure on $(E,\mathcal{E})\otimes (E,\mathcal{E})$ define by
 $$
\mathbf{\mu}(\mathbf{A})=\mu(\{a \in E | (a,a) \in \mathbf{A}\}).
$$
If $\mu_1=\mu_2$, we set
$\P = \mathbf{\mu_1}$.
Else we set
\begin{equation}\label{Prf_norm_var_sens_b}
\P = \mathbf{\mu_1 \wedge \mu_2} + \frac{1}{\norm{\mu_2-\mu_1}_{var}}(\mu_1-\mu_2)^+\otimes (\mu_2- \mu_1)^+.
\end{equation}
Noticing that $a=a\wedge b+(a- b)^+$ and using $\norm{\mu_1-\mu_2}_{var}= (\mu_1 - \mu_2)^+(E)$, we obtain that
$\P(. \times E)=\mu_1 \wedge \mu_2+(\mu_1-\mu_2)^+=\mu_1$ and $\P(E \times .)=\mu_2$. Thus if we denote by
$(Y_i)_i$ the projectors, we obtain that $(Y_i)_i$ is a coupling of $(\mu_i)_i$. Moreover,
$$
\P (Y_1=Y_2, Y_1 \in A)=(\mu_1 \wedge \mu_2)(A).
$$
 So it is the desired maximal coupling
\carre
\begin{Remark}\label{Rq_mesurable}
Moreover, in all this article, we admit that the maximal coupling $(Y_i(u_0^i))_i$  could be chosen such that
 $(Y_i(u_0^1,u_0^2))_i$ depend measurably on the initials conditions $(u_0^i)_i$. The idea is the following.
 Since we only work in nice spaces, we can consider that we are working on the real line. It can be seen that the laws we use
depend measurably on $(u_0^i)_i$ and then the law define by \eqref{Prf_norm_var_sens_b} will do it too. Then its
 repartition function $F_{(u_0^1,u_ 0^2)}$ is measurable too and finally the pseudo-inverse of the repartition function
$F_{(u_0^1,u_ 0^2)}^{-1}$ is measurable with respect to $(u_0^1,u_ 0^2)$. We consider $([0,1],\mathcal{B}_{[0,1]},\lambda)$,
 where $\lambda$ is the Lebesgue measure and we set $Y_i(u_0^1,u_0^2,\omega)=F_{(u_0^1,u_ 0^2)}^{-1}(\omega)$.
 Then $(Y_i)_i$ is measurable with respect to $(u_0^1,u_0^2,\omega)$ and for every $(u_0^1,u_ 0^2)$, it is a coupling
 of $(\mu_i(u_0^i))_i$.  For a proof see \cite{KS}.
\end{Remark}

\section{ Proof of Proposition \ref{Prop_Matt}}
We set
$$
\Omega= E^2,\quad \F=\mathcal{B}(E^2),
$$
and $V_i$ the $i^{th}$ projector on $\Omega$:
$$
V_i(v_1,v_2)=v_i, \quad i=1,2.
$$

Let $(U_1,U_2)$ be a coupling of $(\mu_1,\mu_2)$.

In order to establish Proposition  \ref{Prop_Matt}, we  build a probability measure $Q$
 on $(\Omega,\F)$ such that
\begin{equation}\label{Eq_Matt_a}
\Espace
\left \{
\begin{array}{ll}
\alpha ) & Q(\cdot\times E)=\mu_1, \quad Q(E\times\cdot)=\mu_2,\\
\beta  ) & Q(f_0(V_1)=f_0(V_2))\geq \left(\nu_1\wedge\nu_2 \right)(E).
\end{array}
\right.
\end{equation}

Then $(V_1,V_2)$ seen as a couple of random variables defined on $\left( \Omega, \F,Q \right)$ is a coupling of
$(\mu_1,\mu_2)$ such that $(f(V_1),f(V_2))$ is a maximal coupling of $(\nu_1,\nu_2)$.

Recall that
\begin{equation}\label{Eq_Matt_un}
\nu_i=\nu_1\wedge\nu_2+\left((-1)^i(\nu_1-\nu_2) \right)^+,\quad i= 1,2,
\end{equation}
and that since $E$, $F$ are polish spaces, then there exists a version of \mbox{$\P(U_i\in A\,|\,f_0(U_i)=x)$} which is measurable for
 any $A\in \mathcal B (E)$ and which is probability measure for any $x \in F$. Moreover
\begin{equation}\label{Eq_Matt_deux}
\mu_i(A)= \int_F \P(U_i\in A \,|\,f_0(U_i)=x) \nu_i(dx) ,\quad i= 1,2,
\end{equation}

Combining \eqref{Eq_Matt_un} and \eqref{Eq_Matt_deux}, we obtain
\begin{equation}\label{Eq_Matt_trois}
\mu_i=\mu_i^s + \mu_i^r,\quad i= 1,2,
\end{equation}
where
$$
\Espace
\begin{array}{rclr}
\mu_i^s(A) &=& \int_F \P(U_i\in A \,|\,f_0(U_i)=x) \left(\nu_1\wedge\nu_2 \right)(dx), &i= 1,2,\\
\mu_i^r(A) &=& \int_F \P(U_i\in A \,|\,f_0(U_i)=x) \left((-1)^i(\nu_1-\nu_2) \right)^+(dx) ,& i= 1,2.
\end{array}
$$
Remark that
\begin{equation}\label{Eq_Matt_quatre}
\Espace
\left \{
\begin{array}{l}
\mu_i^s, \mu_i^r \geq 0,\quad i= 1,2,\\
\mu_i^s(E) = \left(\nu_1\wedge\nu_2 \right)(E), \\
\mu_i^r(E) = \norm{\nu_1-\nu_2}_{var}.
\end{array}
\right.
\end{equation}

Taking into account \eqref{Eq_Matt_trois} and \eqref{Eq_Matt_quatre}, we can write problem \eqref{Eq_Matt_a} in the form
\begin{equation}\label{Eq_Matt_b}
\Espace
\left \{
\begin{array}{ll}
\lefteqn{\textrm{Find $r,s$ two positive measures on $(\Omega,\F)$ such that}}\\
i ) & s(\cdot\times E)=\mu_1^s, \quad s(E\times\cdot)=\mu_2^s,\\
ii )  & r(\cdot\times E)=\mu_1^r, \quad r(E\times\cdot)=\mu_2^r,\\
iii) & s(f_0(V_1)\not=f_0(V_2))=0.
\end{array}
\right.
\end{equation}

Once \eqref{Eq_Matt_b} is true, we can set
$$
Q=r+s.
$$
Then \eqref{Eq_Matt_a}$\alpha)$ is an obvious consequence of \eqref{Eq_Matt_trois}. Furthermore, since $r\geq 0$, then
\eqref{Eq_Matt_b}iii), \eqref{Eq_Matt_b}i) and \eqref{Eq_Matt_quatre} gives
 $$
Q(f_0(V_1)=f_0(V_2))\geq s(f_0(V_1)=f_0(V_2))=s(\Omega)=\mu_i^s(E)=\left(\nu_1\wedge\nu_2 \right)(E).
$$

Now we build $r$ by setting
$$
r=\frac{1}{\norm{\nu_1-\nu_2}_{var}} \mu_1^r\times\mu_2^r.
$$
Notice that $r\geq 0$ and  \eqref{Eq_Matt_b}ii) are obvious consequence of \eqref{Eq_Matt_quatre}.

Now we build $s$ by setting
$$
s(A\times B)=\int_F \P(U_1\in A \,|\,f_0(U_1)=x)\times\P(U_2\in B \,|\,f_0(U_2)=x) \left(\nu_1\wedge\nu_2 \right)(dx).
$$
Notice that \eqref{Eq_Matt_b}i) and \eqref{Eq_Matt_b}iii) are obvious.


\end{document}